\documentclass[11pt]{article}
\usepackage{amsmath, amssymb, amsfonts, verbatim, amsthm, xcolor}
\usepackage{graphicx}
\newtheorem{theorem}{Theorem}[section]

\newtheorem{proposition}[theorem]{Proposition}
\newtheorem{corollary}[theorem]{Corollary}
\newtheorem{lemma}[theorem]{Lemma}

\newtheorem{conjecture}[theorem]{Conjecture}
\newtheorem{remark}[theorem]{Remark}

\def\qed{\hfill $\Box$\medskip}
\def\diag{{\rm diag}\,}
\def\span{{\rm span}\,}

\def\ba{{\mathbf a}}
\def\bb{{\mathbf b}}
\def\bc{{\mathbf c}}
\def\bd{{\mathbf d}}
\def\Re{{\rm Re}}

\def\IC{{\mathbb{C}}}
\def\IR{{\mathbb{R}}}

\def\IN{{\mathbb{N}}}

\def\tr{{\rm tr}}

\def\conv{{\rm conv}}
\begin{document}
\openup .92\jot
\title{Numerical ranges  of cyclic shift matrices}
\author{Mao-Ting Chien$^{\rm a}$, Steve Kirkland$^{\rm b}$, Chi-Kwong Li$^{\rm c}$, Hiroshi Nakazato$^{\rm d}$}
\date{}
\maketitle
\begin{itemize}
\item[$~^{\rm a}$]  Department of Mathematics, Soochow University, Taipei 11102, Taiwan (mtchien@scu.edu.tw)
\item[$~^{\rm b}$] Department of Mathematics, University of Manitoba, Winnipeg, R3T 2N2, Manitoba, Canada (Stephen.Kirkland@umanitoba.ca)
\item[$~^{\rm c}$] Department of Mathematics, College of William and Mary, Williamsburg, Virginia 23187-8795, USA (ckli@math.wm.edu)
\item[$~^{\rm d}$] Faculty of Science and Technology, Hirosaki University,  Hirosaki 036-8561, Japan \\ (nakahr@hirosaki-u.ac.jp)
\end{itemize}

\begin{abstract}
We study the numerical range of an $n\times n$ cyclic shift matrix, which can be viewed as 
the adjacency matrix of a directed cycle with $n$ weighted arcs. In particular, we consider 
the change in the numerical range if the weights are rearranged or perturbed.  In addition 
to obtaining some general results on the problem,  a permutation of the given weights is 
identified such that the corresponding matrix  yields the largest 
numerical range (in terms of set inclusion),  
for $n \le 6$. We conjecture that the maximizing pattern 
extends to general $n\times n$ cylic shift matrices.  For $n \le 5$, we also 
determine permutations such that the corresponding cyclic shift matrix yields
the smallest numerical range.
\end{abstract}


{\bf AMS Classification Numbers.} 15A60, 05C50, 05C38.

\section{Introduction}

Let $M_n$ be the set of $n\times n$ complex matrices. The numerical range of $A\in M_n$ is defined and denoted by
$$W(A)= \{x^*Ax: x \in \IC^n, x^*x = 1\},$$
 and is a useful tool for studying matrices; for example, see \cite[Chaper 1]{HJ}.
In particular, there is an interesting interplay between the geometric properties of 
$W(A)$ and the algebraic properties of $A$. The T\"{o}plitz-Hausdorff theorem  
asserts that $W(A)$ is a convex set. It is known that $W(A) = \{\lambda\}$
if and only if $A = \lambda I$; $W(A) \subseteq \IR$ if and only if  
$A$ is Hermitian; $W(A)$ is a polygon (with interior) if and only if
$A$ is unitarily similar to $D \oplus B$ for a diagonal matrix $D\in M_k$ (for some $k=1, \ldots, n$) such that
$W(B) \subseteq W(A) = W(D) = \conv ( {\rm{eig}}(D))$, the convex hull of spectrum of $D$; e.g, 
see \cite{HJ}.

In this paper, we study the numerical range of a  {\it cyclic shift matrix}
defined as follows. Let $a_1, \dots, a_n\in \IC$, and let 
$$A = S(a_1, \dots, a_n) = \sum_{j=1}^{n-1} a_j E_{j,j+1}+ a_n E_{n,1}$$
be the corresponding cyclic shift matrix. Here 
$\{E_{1,1}, E_{1,2}, \dots, E_{n,n}\}$ is the standard basis for $M_n$.
There are a number of interesting papers on properties of the numerical ranges  of 
cyclic shift matrices (see \cite{CN1, CN5, GTW, GW, Tsai-1, Tsai-2, TGW, WW}). 
In general, if $A \in M_n$ is such that each of its 
 rows and  columns  has at most one nonzero entry,
then $A$ is permutationally similar to  $A_1 \oplus \cdots \oplus A_r$ (possibly with an additional direct summand that is a zero matrix), 
where each $A_j$ is a cyclic shift matrix. Hence,  our results can be 
used to describe $W(A)$ in terms of $ \conv (\cup_{j=1}^r W(A_j))$.

Now, suppose $A = S(a_1 \dots, a_n)$. 
There is a diagonal unitary matrix $U$ such that for some $t \in \IR,$ 
$U^*AU = e^{it} \hat A$, where $\hat A = \sum_{j=1}^{n-1} |a_j| E_{j,j+1} + |a_n| E_{n,1},$ 
so that $W(A) = e^{it} W(\hat A)$. 
Hence, we may focus our attention on the case that  $a_1, \dots, a_n \ge 0$.
By the result in \cite{LT}, if $a_1, \dots, a_n > 0$, we have
$$W(A) = e^{i2\pi/n} W(A).$$ On the other hand, 
if $a_j = 0$ for some $j$, then $W(A)$ is a circular disk centered at the origin with 
radius equal to the largest eigenvalue of $\Re(A) = (A+A^t)/2$.
In fact, if $k$ of the values $a_1, \dots, a_n$ are equal to  zero, then
up to a permutation similarity, $A$ has the form  
$A_1\oplus \cdots \oplus A_k$ where each $A_j$ has positive entries on the superdiagonal, 
and $W(A_j)$ is a circular disk centered at the origin 
with radius equal to the largest eigenvalue, i.e., the spectral radius,  
of $(A_j + A_j^t)/2$. If the largest eigenvalue of $A+A^t$ corresponds to the 
direct summand $A_1+A_1^t$, then $W(A) = W(A_1)$.

Let $A = S(a_1, \dots, a_n)$ and $B = S(b_1, \dots, b_n)$
for $a_j, b_j \ge 0$, $j = 1, \dots, n$. We are interested in comparing 
$W(B)$ and $W(A)$. In particular, if
$(b_1, \dots, b_n)$ is obtained from $(a_1, \dots, a_n)$ by a 
permutation $\sigma$ of the entries, we write $B = A_\sigma$. 

Let $r(A) = \max\{|\mu|: \mu \in W(A)\}$ be the numerical radius of 
$A \in M_n$. 
In \cite{CW}, it is shown that if $A = S(a_1, a_3, a_5, \dots, a_6, a_4, a_2)$
with $0 \le a_1 \le \cdots \le a_n$,  then
$$r(A_\sigma) \le r(A) \qquad \hbox{ for any permutation } \ \sigma.$$
In particular, if $a_1 = 0$, then $W(A_\sigma)$ is a circular disk, and 
\begin{equation}\label{conjecture}
W(A_\sigma) \subseteq W(A) \qquad \hbox{ for any permutation } \ \sigma.
\end{equation}
We formulate the following conjecture. 
	
\begin{conjecture}\label{conj1} 
	Suppose that $n \ge 2,$  that  $0 \le a_1 \le \cdots \le a_n$, and let $A = S(a_1, a_3, a_5, \dots, a_6, a_4, a_2).$ Then {\rm (\ref{conjecture})} holds
	for any $\sigma \in S_n$.
\end{conjecture}

Conjecture \ref{conj1} holds trivially
for $n=2,3$ and in this paper, we confirm the conjecture for  $n=4,5,6$. 
We also study permutations $\tau$ such that
$W(A_\tau) \subseteq W(A_\sigma)$ for all $\sigma$. 
For $n = 4$, we have $\tau = (1,3,2,4)$;
for $n = 5$, we have $\tau = (1,4,3,2,5)$.

Note that  $S(a_1, a_2, \dots, a_n)$ is permutationally similar to 
$S(a_n,a_1, a_2,\dots, a_{n-1}),$ which is obtained by a clockwise permutation of the weights 
(cf. \cite{Tsai-2}). 
More generally,
letting $C = S(1,\dots, 1)$, 
we may replace $A_\sigma$ by $C^r A_\sigma C^{-r}$ (for some suitable choice of $r$) and 
assume that the $(1,2)$ entry of $A$ is the minimum
of the values $a_1, \dots, a_n$.
In our discussion, we often  assume that $A = S(a_1, \dots, a_n)$ with 
$a_1,  \cdots, a_n > 0$.  Letting 
$a_j \rightarrow 0$, may then recover  the corresponding results  when some $a_j = 0$. 

The rest of the paper is organized as follows.
In Section 2, we  present some preliminary results. In Section 3, we consider 
extreme matrices $A = S(a_1, \dots, a_n)$ attaining the maximum and minimum
numerical ranges under some restrictions on $a_1, \dots, a_n$.
In Sections 4--6, we compare $W(A)$ and $W(B)$ for 
$A = S(a_1,\dots, a_n)$ and $B = S(b_1, \dots, b_n)$ with emphasis on the
cases that $(b_1, \dots, b_n)$ is a permutation of $(a_1, \dots, a_n)$.
We mention some related results and future research directions in Section 7.

\section{Preliminaries}

Denote by $\lambda_1(H)\ge \cdots \ge \lambda_n(H)$
the eigenvalues of a Hermitian matrix $H \in M_n$.
For $\theta \in [0, 2\pi)$,  $A \in M_n$, consider  
$ \Re(e^{i\theta }A) = (e^{i\theta}A + e^{-i\theta}A^*)/2$. 
We have the following; see \cite[Theorem 1.5.12]{HJ}.

\begin{proposition}
\label{maxeig}  Let $A \in M_n$. Then
$$W(A) = \cap_{\theta \in [0, 2\pi)} 
\{ \mu \in\IC: \Re(e^{i\theta} \mu) \le \lambda_1(\Re(e^{i\theta}A))\}.$$
As a result, for any $B \in M_n$,
$W(B) \subseteq W(A)$ if and only if $\lambda_1(\Re(e^{i\theta} B)) \le
\lambda_1(\Re(e^{i\theta} A))$ for all $\theta \in [0, 2\pi).$
\end{proposition}

\begin{proposition} 
\label{observations}
Let $A = S(a_1, \dots, a_n)$ be such that  $a_j >0$ ($j=1,2, \ldots, n$), and fix $\theta \in [0, 2\pi).$  
Then each eigenvalue of $\Re(e^{i\theta }A)$ has multiplicity at most 2.

\begin{itemize}
\item[{\rm (a)}] Suppose $n=2k$ is even. Then $\Re(e^{i\theta }A) $ 
is permutationally similar to a matrix of the form 
$\begin{pmatrix}0_k & R\cr  R^* & 0_k\cr\end{pmatrix}$, where $0_k$ is the $k \times k$ zero matrix.  The eigenvalues of 
$\Re(e^{i\theta }A) $ have the form $\pm s_j$, where $s_1 \ge \cdots \ge s_k\ge 0$ are the 
singular values of $R$. 
 If  $s_j = s_{j+1}$
for some $j$ then $-s_j = -s_{j+1}$ is another repeated eigenvalue.
Moreover,
if   $\det(zI - (iA-iA^*)) =  \prod_{j=1}^k(z^2-\alpha_j^2)$, then 
  $$\det(zI-2\Re(e^{i\theta }A))
=  \prod_{j=1}^k(z^2-\alpha_j^2)+(-1)^k 2\prod_{j=1}^n a_j- 2\Big(\prod_{j=1}^n a_j\Big) \cos(n\theta).$$

\item[{\rm (b)}] Suppose $n = 2k+1$ is odd. The matrix obtained from 
$\Re(e^{i\theta }A) $ by deleting the last row and last column has 
$k$ positive and $k$ negative eigenvalues. The matrix
$\Re(e^{i\theta }A)$ has at least $k$ positive and $k$ negative eigenvalues.
Moreover, if 
$\det(zI - (iA-iA^*)) = z \prod_{j=1}^k(z^2-\alpha_j^2)$, then 
$$\det(zI-2\Re(e^{i\theta }A)) =
z \prod_{j=1}^k(z^2-\alpha_j^2) -2\Big( \prod_{j=1}^n a_j \Big) \cos(n\theta).$$
\end{itemize}
\end{proposition}

{\it Proof.} 
Removing the first row and the last row of the matrix $\lambda I - \Re(e^{i\theta }A)$ 
 yields an $(n-2)\times n$ matrix in row echelon form
with positive $(j,j)$ entry for $j = 1, \dots, n-2$. So, 
the kernel of $\lambda I - \Re(e^{i\theta }A)$ has dimension at most 2.

(a)  If we apply a permutation similarity to $\Re(e^{i\theta }A)$ and put the rows and columns
with odd indices as the first $k$ rows and columns,
we see that $\Re(e^{i\theta }A)$ has the asserted block form.   The assertions about 
the eigenvalues follow. 

To prove the formula for $\det(zI-2 \Re(e^{i\theta }A))$,
let $D = \diag(1, e^{i\theta}, e^{2i\theta}, \dots, e^{i(n-1)\theta})$.
Then $D(2\Re(e^{i\theta }A))D^*$ is the same as $A+A^t$ except that the $(1,n)$
and $(n,1)$ entries are equal to $a_n e^{-in\theta}$ and $a_n e^{in\theta}$.
Now, $\det(zI - 2\Re(e^{i\theta }A)) = z^n+ \sum_{j=1}^n \gamma_j(\theta) z^{n-j}$
such that $\gamma_j(\theta)$ is the sum of the $j\times j$ principal minors
of $2\Re(e^{i\theta }A)$. One may readily check  that  
$\gamma_1(\theta), \dots, \gamma_{n-1}(\theta)$
are independent of $\theta$, and $\gamma_n(\theta) 
= \gamma_n(\pi/2)  +(-1)^k 2a_1 \cdots a_n -2 a_1 \cdots a_n \cos(n\theta)$.
The result follows.

(b) Removing the last row and the last column of $\Re(e^{i\theta }A)$,
we get a tridiagonal Hermitian matrix that is permutationally similar to 
a matrix of the form  $\begin{pmatrix} 0 & R \cr R^* & 0 \cr\end{pmatrix}$.
We get the conclusion on the eigenvalues of the matrix.
We can use the arguments in the proof of (a) to conclude that 
$\det(zI - 2\Re(e^{i\theta }A)) = \det(I -(iA-iA^*)) - 2a_1 \cdots a_n\cos(n\theta)$,
and $\det((iA-iA^*) ) = 0$. The conclusion follows.
\qed

It is interesting to note that for $C = S(1, \dots, 1)$, then 
$C+C^t$ has eigenvalues $2\cos (2k\pi/n)$ for $k = 1, \dots, n$. 
Thus, if $n \ge 3$ is odd, there are  $(n-1)/2$ distinct  
eigenvalues  with multiplicity 2, while if $n \ge 4$ is even, there are  $(n-2)/2$ distinct eigenvalues with multiplicity 2.  
On the other hand, if $A = S(a_1, \dots, a_{n-1},a_n)$ with 
$a_1 \cdots a_{n-1} > 0 = a_n$, then $A+A^t$ necessarily has $n$ distinct 
eigenvalues.

The following
proposition, which is of independent interest (especially in 
the graph theory context) shows that 
 one can change 
$(a_{n-1}, a_n)$ to produce a nonnegative matrix $\hat A$ 
so that the smallest eigenvalue of $\hat A + \hat A^t$ 
has multiplicity 2.

\begin{proposition} Suppose $n$ is odd, and 
$A = S(a_1, \dots a_{n-1}, a_n)$ with $a_1 \cdots a_{n-1} > 0 = a_n$.
Then there are $\hat a_{n-1}, \hat a_n > 0$ such that
for $\hat A = S(a_1, \dots, a_{n-2}, \hat a_{n-1}, \hat a_n)$,
the smallest eigenvalue of $\hat A + \hat A^t$
has multiplicity 2. 
\end{proposition}

{\it  Proof.}  Let $n = 2k+1$, and $A = S(a_1,\dots, a_{n-1},a_n)$
with $a_1 \cdots a_{n-1} > 0 = a_n$. Then removing the last row and the
last column of the matrix $A+A^t$, we get $B$ with $n-1$ distinct eigenvalues.
Let $v\in \IR^{n-1}$ be an eigenvector of $B$ corresponding to the smallest eigenvalue
$\mu$. Then the first entry $v_1$ and last entry $v_{n-1}$ of $v$
have opposite signs. (This follows from Perron--Frobenius theory: the graph associated with $B$ is bipartite, so $-\mu$ is the spectral radius of $B$; the conclusion that $v_1v_{n-1}<0$  now follows from the fact that vertices $1$ and $n-1$ are in different subsets of the bipartition.)
Let $b = (|v_{n-1}|, 0,\dots,0, |v_1|)^t\in \IR^{n-1}$, 
and  form the matrix
$M(x)=\begin{pmatrix}
B  & xb\cr xb^t &  0\cr\end{pmatrix}$, 
with $x\in \IR$. Then $\mu$ 
is an eigenvalue of $M(x)$ for any $x \in \IR$, 
with corresponding eigenvector $\begin{pmatrix} v \cr 0\cr\end{pmatrix}$.
For $x=0, M(0)$ 
has same eigenvalues as $B$, with an extra 0 appended. 
As $x\rightarrow\infty$ 
one of the eigenvalues of $M(x)$ goes to 
$-\infty$. By the Intermediate Value Theorem, 
there must be a value $x_0>0$ such that $M(x_0)$ has 
$\mu$ as a double eigenvalue. 
\qed

In connection to the comparison of the maximal eigenvalue of Proposition \ref{maxeig}, 
the following simple lemma is useful in our study.

\begin{lemma} \label{poly}
Suppose that $f(t), g(t)$ are monic polynomials of degree $n$ and that  each
has $n$ real zeros. 
Suppose further that  $f(t) = g(t)-h(t)$ for some polynomial $h(t)$ of degree
smaller than $n$ such that  one of the following conditions holds:

\begin{itemize}
\item[{\rm (a)}] 
$h(t_1) \ge 0$, where $t_1$ is the largest zero of $g(t)$,
\item[{\rm (b)}] 
$h(t) >0$ for $t \geq t_2$, where $t_2$ is the largest zero of $f(t)$. 
\end{itemize}
Then the  largest positive zero of $f(t)$ is larger than or equal to  that of $g(t)$.
\end{lemma}

{\it Proof.}
(a) Note that $f(t_1) \le 0$. Thus, the largest zero of $f(t)$ is larger than or
equal to $t_1$. 

(b) If $h(t) >0$ for $t \geq t_2$, then $g(t)=f(t) +h(t) > 0$ for $t \geq t_2$ and hence $t_1 < t_2$. \qed

\begin{proposition}\label{alpha-beta}  
Suppose $A = S(a_1, \dots, a_n)$ and 
$B = S(b_1, \dots, b_n)$ with $a_j, b_j > 0$ for $j = 1, \dots, n$.
Let $\alpha =  \prod_{j=1}^n a_j$, $\beta = \prod_{j=1}^n b_j$,
There are monic polynomials $f(z), g(z)$ of degree $n$
independent of $\theta \in [0,2\pi)$ such that  for any such $\theta,$ 
$$\det(zI -2 \Re(e^{i\theta}A)) = f(z)-2 \alpha \cos n \theta \quad \hbox{ and }
\quad
\det(zI -2 \Re(e^{i\theta}B)) = g(z)-2 \beta  \cos n\theta.$$
\medskip\noindent
Let $z_t$ 
be the largest zeros of $f(z)-2\alpha t$ for $t \in [-1,1]$ and the largest zero of 
$g(z) + 2\beta$ is not larger than  $z_{-1}$. Then $W(B)\subseteq W(A)$ if and only if
$g(z)/\beta \ge f(z)/\alpha$ for $z \in [z_{-1}, z_1]$. 
\end{proposition}

{\it Proof.}
Note that for any $n \times n$ matrix $M$, $\det(zI_n - M) = z^n + \sum_{j=1}^n m_j z^{n-j}$
 where  
$(-1)^j m_j$ is the sum of the $j\times j$ principal minors of $M$.
One readily checks that the coefficients of $z^{n-j}$ of the polynomial $\det(zI -2\Re(e^{i\theta}A))$ for $2 \leq j \leq n-1$
are independent of $\theta$ and correspond to the (same) sum of principal minors.
The constant term is given by 
  $$\det(2 \Re(e^{i\theta}A)) = 2 \alpha \cos(n\theta)$$ 
for odd $n,$  and
  $$\det(2 \Re(e^{i\theta}A)) =(-1)^{n/2}(a_1^2 a_3^2 \cdots a_{n-1}^2 +a_2^2 a_4^2 \cdots a_n^2) -2 \alpha \cos(n \theta)$$
for even $n$. So we conclude that $\det(zI-2\Re(e^{i\theta}A))$ and 
$\det(zI-2\Re(e^{i\theta}B))$ have the asserted forms.

 Let $t = \cos n\theta \in [-1,1]$. 
Then $f(x) - 2\alpha t$ has $n$ real zeros for every $t \in [-1,1]$, 
and $f(x)$ is a convex increasing function for $x \in [z_{-1},\infty)$.
Similarly,  $g(x) - 2\beta t$  has $n$ real zeros for every $t \in [-1,1]$, 
and $f(x)$ is a convex increasing function for $x \in [\hat z_{-1},\infty)$,
where $\hat z_{-1}$ is the largest zero of $g(x) + 2\beta$.
By our assumption, $\hat z_{-1} \le z_{-1}$.

Note that $W(B) \subseteq W(A)$ if and only if 
for every $t \in [-1,1]$, the 
largest root of $f(x) - 2\alpha t = 0$ always lies on the right side of that of
$g(x) - 2\beta t = 0$.  Equivalently, in $[z_{-1},z_1]$, the curve $g(x)/\beta$ is always above
that of $f(x)/\alpha$, and the result follows. \qed

 Note that $f(z)$ has the form $(x^2-c_1^2)\cdots (x^2-c_k^2)$
or $x(x^2-c_1^2)\cdots (x^2-c_k^2)$ depending on whether $n = 2k$ or $n=2k+1$.

\section{Cyclic shift matrices with extreme numerical ranges}

In the following, we show that some weighted cycles and paths will 
yield the smallest numerical ranges under certain restrictions on the weights.

\begin{theorem} \label{3.1} Suppose $A = S(a_1, \dots, a_n)$ such that
$a_1\cdots a_n = 1$ and consider $C = S(1,\dots, 1)$.
Then $W(C) \subseteq W(A)$. Equality holds if and only if
$|a_1| = \cdots = |a_n| = 1$.
\end{theorem}

{\it Proof.}  Note that $\det(\lambda I-A) = \lambda^n -1$.
Thus, $A$ has eigenvalues  $e^{i2k\pi/n}$ for $k = 1, \dots, n$.
By convexity of $W(A)$ and the fact that $W(A)$ contains all the 
eigenvalues of $A$, we see that 
$W(C) = \conv\{e^{i2k\pi/n}:  1\le k \le n\} \subseteq W(A)$.
If  $W(A) = W(C)$, then all eigenvalues of $A$ are non-differentiable
boundary  points of $W(A)$; e.g., see \cite{HJ}. 
Hence, $A$ is normal, and equals $C$.
\qed

\medskip
Below is  a similar result for 
$A = S(a_1, \dots, a_n)$ with $a_n = 0$. In this case, $W(A)$ is always a 
circular disk centered at the origin with radius $r(A)$.
We will determine the smallest radius if $\prod_{j=1}^{n-1} a_j = 1$.

\begin{theorem} \label{3.2}  
Suppose that $A = S(a_1, \dots, a_n)$ for some 
$a_1, \dots, a_n \in \IC$ with 
$|\prod_{j=1}^{n-1} a_j| = 1$ and $a_n = 0$.
If $|A| = S(|a_1|, \dots, |a_{n-1}|, 0)$, then $W(A) = W(|A|)$ is a circular disk
centered at the origin with radius $\lambda_1(|A|+|A|^t)/2$ 
such that
$$\lambda_1(A+A^*) = \lambda_1(|A|+|A|^t)  \ge 2^{(n-2)/(n-1)}.$$ 
The equality holds if and only if
$|a_1|=|a_{n-1}|= 2^{(n-3)/(2n-2)}$, $|a_j|=2^{-1/(n-1)}$ for  $j=2, ... , n-2$,
so that 
 $x = (1/\sqrt{2n-2}, 1/\sqrt{n-1}, 1/\sqrt{n-1},\ldots,1/\sqrt{n-1}, 1/\sqrt{2n-2})^t$ is the corresponding
(Perron) eigenvector for $|A| + |A|^t$.
\end{theorem}

{\it Proof.} 
We know that $W(A) = W(|A|)$ is a circular disk centered at the origin with 
radius $\lambda_1(|A|+|A|^t)/2$.  So, 
suppose that $A = S(a_1, \dots, a_n)$
gives rise to the minimum $\lambda_1(A+A^t) = r$ with 
$a_1, \dots, a_{n-1} > 0 = a_n$. 
 Consider the function 
  $$f(a_1,\ldots, a_{n-1}, x_1, \ldots, x_n)=\sum_{j=1}^{n-1} a_j x_j x_{j+1}.$$
 The  minimizing matrix $A = S(a_1, \dots, a_n)$ can be obtained by finding the minimax value 
$$\min_{a_1,\ldots, a_{n-1}}\max_{x_1,\ldots,x_n} f(a_1,\ldots, a_{n-1}, x_1, \ldots, x_n)$$
over the conditions
  $a_1 a_2 \cdots a_{n-1}=1$ and   $\sum_{j=1}^n x_j^2 =1$,
and $x_j \geq 0$, $1 \leq j \leq n$, $a_j \geq 0$, $1 \leq j \leq n-1$. 
By substituting $a_{n-1}=1/(a_1a_2\cdots a_{n-2})$ and $x_{n}=\sqrt{1-(x_1^2+\cdots +x_{n-1}^2)}$ into the function $f(a_1,\ldots, a_{n-1}, x_1, \ldots, x_n)$, the minimax value of the function $f(a_1,\ldots, a_{n-1}, x_1, \ldots, x_n)$ reduces to the  minimax value of the function 
$$F(a_1,\ldots, a_{n-2}, x_1, \ldots, x_{n-1})=\Big(\sum_{j=1}^{n-2} a_j x_j x_{j+1}\Big)+{1\over a_1a_2\cdots a_{n-2}} x_{n-1}\sqrt{1-(x_1^2+\cdots+x_{n-1}^2)},$$
i.e.,
$$\min_{a_1,\ldots, a_{n-2}}\max_{x_1,\ldots,x_{n-1}} F(a_1,\ldots, a_{n-2}, x_1, \ldots, x_{n-1})$$
 over the variables $a_1,a_2,\ldots, a_{n-2}>0$ and $x_1^2+x_2^2+\cdots+x_{n-1}^2<1$.    

Since changing $a_i$ will increase $F(a_1,\ldots, a_{n-2}, x_1, \ldots, x_{n-1})$, and  changing $x_i$ will decrease $F$.
It follows that at the minimum,  all the partial derivatives of  $F(a_1,\ldots, a_{n-2}, x_1, \ldots, x_{n-1})$ will be 0. 
Then we see that
$$0 = x_jx_{j+1} -{1\over a_j} \alpha,\ \quad j = 1, \dots, n-2,\eqno(3.1)$$
$$0 = a_1 x_2 -{x_1\over \beta} = a_{n-2}\ x_{n-2} -  {x_{n-1}\over \beta}+{\alpha\over x_{n-1}}= a_{j-1} x_{j-1} + a_j x_{j+1}-{x_j\over \beta},     \quad j = 2, \dots, n-2.\eqno(3.2)$$
where $\alpha={x_{n-1}\over  a_1a_2\cdots a_{n-2}} \sqrt{1-(x_1^2+\cdots +x_{n-1}^2)}$
and ${1\over\beta}={x_{n-1}\over (a_1\cdots a_{n-2}) \sqrt{1-(x_1^2+\cdots +x_{n-1}^2)}}$

By the equation (3.1), 
$$a_j = \alpha/(x_jx_{j+1}), \quad j = 1, \dots, n-2. \eqno(3.3)$$
Substituting (3.3) into the  equation (3.2), we get
$$\alpha\beta = x_1^2 =x_n^2=x_j^2/2, \quad j = 2, \dots, n-1.$$
Consequently, the condition $\sum_{j=1}^n x_j^2=1$ yields $\alpha\beta=1/(2n-2)$, and thus
$x_1=x_n=1/\sqrt{2n-2}$ and $x_j=1/\sqrt{n-1}$ for $j=2,3,\ldots, n-1$.
Now, 
 $${1\over\beta}={x_{n-1}\over (a_1\cdots a_{n-2}) \sqrt{1-(x_1^2+\cdots +x_{n-1}^2)}}={x_{n-1}\over (a_1\cdots a_{n-2}) x_n}=a_{n-1}\sqrt{2}.$$
Thus $a_{n-1}={1\over \sqrt{2}\beta}$. Then, together with $a_1,\ldots, a_{n-2}$ in (3.3), we bave that
$$1=a_1\cdots a_{n-2}a_{n-1}={\alpha\over x_1x_2}{\alpha\over x_2x_3}\cdots {\alpha\over x_{n-2}x_{n-1}}{1\over \sqrt{2}\beta}$$
$$= {\alpha\over x_1x_2}{\alpha\over x_2x_3}\cdots {\alpha\over x_{n-2}x_{n-1}}{\alpha\over \sqrt{2}\alpha\beta}$$
$$={\alpha^{n-1}\over {1\over \sqrt{2}(n-1)^{n-2}}}{2n-2\over \sqrt{2}}=\alpha^{n-1} 2 (n-1)^{n-1}.$$
Hence
$$\alpha={1\over n-1}2^{-{1\over n-1}}.$$
Again by (3.3), we have
$$a_1=2^{(n-3)/(2n-2)}, \ {\rm and} \ a_j=2^{-1/(n-1)}, j=2,3, \ldots, n-2.$$
From the identity
$a_{n-1}={1\over \sqrt{2}\beta}$, we obtain that
$$a_{n-1}={\alpha\over \sqrt{2} \alpha\beta}={2n-2\over \sqrt{2}}{1\over n-1}2^{-{1\over n-1}}=2^{(n-3)/(2n-2)}.$$ 
In this case, $f(a_1,\ldots, a_{n-1}, x_1, \ldots, x_n)=2^{-1/(n-1)}$. 
Since $f(k^{n-2},1/k,\ldots,1/k,1/\sqrt{n},\ldots,1\sqrt{n})>k^{n-2}/n \rightarrow \infty$ as $k\rightarrow \infty$, it follows that $2^{-1/(n-1)}$ is the minimum of $f$. 
 \qed

One may wonder whether we can find $A = S(a_1, \dots, a_n)$ with 
$\prod_{j=1}^n a_j = 1$ with the largest numerical range.
This is not possible as one can let  $A = S(r^{n-1}, r^{-1}, r^{-1}, \dots, r^{-1})$
for $r > 0$ so that $W(A)$ will contain $W(B)$ with $B = r^{n-1} E_{1,2} \in M_2$ such that $W(B)$ is a circular disc centered at the origin with radius $r^{n-1}/2$.

One may consider imposing other conditions on the weights 
of $A = S(a_1, \dots, a_n)$ and determining the matrix 
with extreme numerical range. For instance, suppose 
$A = S(a_1, \dots, a_n)$ with 
Frobenius norm 1, i.e.,  $\sum_{j=1}^n |a_j|^2 = 1$. We have the following.

\begin{theorem}  Suppose $n \ge 3$
and $B = S(b_1, \dots, b_n)$ with $b_1, \dots, b_n \ge 0$ such that
$b_1^2 + \cdots + b_n^2 = 1$.
 If $\prod_{j=1}^n b_j = 0$, then $W(A_0) \subseteq W(B)$, where
$$A_0 = \begin{cases} \sqrt{2/n} S(1,0,\dots,1,0)& 
\hbox{ if } n \hbox{ is even},\cr
\sqrt{2/(n-1)}S(1,0,\dots,1,0, \cos(\theta), \sin (\theta), 0) 
\hbox{ for some } \theta \in [0, \pi/2] 
& \hbox{ if } n \hbox{ is odd}.\cr\end{cases}$$
The equality holds
if and only if
$B$ is permutationally similar to $A_0$. 
\end{theorem}
 
{\it Proof.}
 Suppose $\prod_{j=1}^n b_j = 0$. 
Then we may assume that $b_j \ge 0$ with $b_n = 0$ so 
that $B$ can only have nonzero entries at the $(1,2), \dots, (n-1,n)$ positions.
So, $W(B)$ is a circular disk centered at the origin with radius $r(B)$.
Suppose $B$ is the matrix attaining the minimum numerical radius.

Case 1.  Assume that $n$ is even. By the result in \cite[Theorem 3 (iii)]{JL}, a matrix $X \in M_n$ with Frobenius
norm $\|X\|_F = \tr(XX^*) = 1$ attaining the minimum numerical radius 
must be unitarily similar to $A_0$. We get the desired conclusion.

Case 2. Assume that $n$ is odd. 
We prove the assertion by induction on $k$ if $n=2k+1$.
If $n = 3$, then $\lambda_1(B+B^t) =
\sqrt{b_1^2 + b_2^2} = 1$ so that $W(B)$ is a circular disk centered at the origin with radius $1/2$.
Suppose the assertion holds for $n = 2k-1 \ge 3$.
Consider $n = 2k+1$.
Let $B = S(b_1, \dots, b_{2k}, 0)$
attain the smallest numerical radius, where
$b_1, \dots, b_{2k} \ge 0$ 
satisfying $\sum_{j=1}^{2k} b_j^2 = 1$. 
We may further assume that $b_1 > 0$. Otherwise, replace
$B$ by $PBP^t$ with $P = S(1,\dots, 1)^r$ for a suitable $r \in \IN$.

Note that $W(B) \subseteq W(A_0)$, and  
$A_0$ is a direct sum of $k-1$ matrices
of the form $\gamma S(1,0)\in M_2$
and a matrix of the form $\gamma S(\cos(\theta), \sin(\theta), 0) \in M_3$
with $\gamma = \sqrt{1/k}$ so that
$W(A_0)$ is a circular disk centered at the origin with 
radius $\gamma/2 = 1/\sqrt{4k}$.

Suppose $b_1^2+b_2^2 = \beta^2 > 1/k$. If $C_1$ is the
leading $3\times 3$ submatrix of $B$, then $W(C_1)$ is a circular
disk centered at the origin with radius $\beta/2 > 1/\sqrt{4k}$
so that $r(B) \ge r(C_1) > r(A_0)$,
contradicting the fact that $r(B)$ attains the minimum.
If $b_1^2 + b_2^2 < 1/k$, then $b_3^2 + \cdots + b_{n-1}^2 > 1-1/k$.
If $A_2$ and $C_2$ are the trailing $(n-2)\times (n-2)$ submatrices 
of   $A_0$ and $B$, then $\tr(C_2 C_2^*) \ge \tr(A_2A_2^*)$.
By the induction assumption, $r(A_2)$ is minimum among the matrices 
$C = S(c_1, \dots, c_{n-3}, 0)$  such that 
$C$ and $A_2$ have the same Frobenius norm.
As a result, 
$r(B) \ge r(C_2) > r(A_2) = r(A_0)$, contradicting the fact that $r(B)$
is the minimum. 

Thus,  $b_1^2+b_2^2 = 1/k$, and $C_2$ has the smallest numerical range
among all matrices $S(c_1, \dots, c_{n-3},0)$ in $M_{n-2}$ with the same Frobenius norm as $C_2$. 

If $b_2 = 0, b_1 = \sqrt{1/k}$ so that $B = S(b_1,0) \oplus C_2$ 
since $r(C_2)$ attains the minimum among  
for cyclic shift matrices $S(d_1, \dots, d_{n-3},0) \in M_{n-2}$
with the same Frobenius norm as $C_2$. 
By the induction assumption $C_2$ is permutationally similar to $A_2$.

If $b_2 \ne 0$, then we claim that $b_3 = 0$. Otherwise, we can let 
$B_3 = S(b_1,b_2,0)$ and $B_4 = S(b_1,b_2,b_3,0)$
be the leading $3\times 3$ and $4\times 4$ submatrices of $B$.
We have
 $$\lambda_1(B_4 + B_4^t) > \lambda_1(B_3+B_3^t) = 1/\sqrt{k},$$
contradicting the fact that $r(B_4) \le r(B) \le r(A_0) = 1/\sqrt{4k}$.
Under the assumption, we have $b_2 > 0$. Recall that we assume that $b_1 > 0$.
So, $B = S(b_1, b_2,0) \oplus B_2$, 
where $B_2 = S(b_4, \dots, b_{n-1},0)$ attains the minimum 
numerical radius among matrices with the same Frobenius norm as itself.
Hence, $B_2$ is a multiple of $S(1,0,\dots, 1,0)\in M_{n-3}$ by the result
in the even case.
Thus, the assertion holds for $B \in M_{2k+1}$. \qed

\section{Results for $n = 2,3,4$}

Next, we present some general results when $n = 2,3,4$.

\begin{theorem} Suppose that 
  $A = S(a_1, a_2), B = S(b_1, b_2) \in M_2$ 
 with $a_1, a_2, b_1, b_2 \ge 0$. Then 
$W(B) \subseteq W(A)$ if and only if $|b_1-b_2| \le |a_1-a_2|$ 
 and $b_1+b_2 \le a_1 + a_2$.
 \end{theorem}
 
 {\it Proof.} Note that 
 $W(A)$ is an elliptical disk with major axis $[-(a_1+a_2), (a_1+a_2)]/2$
 and minor axis $i[-|a_1-a_2|,|a_1-a_2|]/2$, and that $W(B)$ has a similar description. 
 The result follows.
 \qed

In the context of Theorem 4.1, note that if  $a_1 > a_2 = 0$, then  $W(B) \subseteq W(A)$ if and only if $b_1+b_2 \le a_1$. 
For $n = 3$, we have the following.
\begin{theorem}
Suppose $A = S(a_1, a_2, a_3), B= S(b_1,b_2,b_3) \in M_3$
 with $a_1, a_2, a_3, b_1, b_2, b_3 > 0$.
Let $\alpha_1 = a_1^2 + a_2^2 + a_3^2, \alpha_0 = a_1 a_2 a_3$,
$\beta_1 = b_1^2 +b_2^2 + b_3^2, \beta_0 = b_1 b_2 b_3$.
 Then 
$W(B) \subseteq W(A)$ if and only if 
for $\xi \in \{1,-1\}$,
the largest zero of 
$z^3 - (a_1^2+a_2^2+a_3^2)z - 2a_1a_2a_3 \xi$
is larger than  or equal to  that of 
$z^3 - (b_1^2+b_2^2+b_3^2)z - 2b_1b_2b_3 \xi$.
Consequently, we have the following.  
\begin{itemize}
\item Suppose $a_1a_2a_3= b_1b_2b_3$. Then $W(B) \subseteq W(A)$ if and only if 
$a_1^2+a_2^2+a_3^2 \ge b_1^2 + b_2^2 + b_3^2$.
\item Suppose $a_1^2+a_2^2 + a_3^2 = b_1^2+b_2^2+b_3^2$.
Then $W(B) \subseteq W(A)$ if and only if $a_1a_2a_3 \ge b_1b_2b_3$.
\item 
For any permutation $\sigma$ of $(1,2,3)$, $W(A) = W(A_\sigma)$. 
\end{itemize}
\end{theorem}

{\it Proof.} 
Note that
$\det(zI-2\Re(e^{i\theta}A))$ and 
$\det(zI-2\Re(e^{i\theta}B))$ have the respective forms:
$$z^3  - \alpha_1 z  - 2\alpha_0\cos(3\theta)
\quad\hbox{ 
and } \quad 
z^3  - \beta_1 z  - 2\beta_0\cos(3\theta).
$$
 Suppose that  $W(B) \subseteq W(A)$.
Then for any $\theta \in [0, 2\pi)$, the largest zero of 
$\det(zI-2\Re(e^{i\theta}A))$ is not less than that of 
$\det(zI-2\Re(e^{i\theta}B))$. In particular, this is true for
$\cos\theta =\{-1,1\}$. 
Conversely, assume that the condition for the largest zero holds for 
$\cos\theta =\{-1,1\}$. If there is $u \in (-1,1)$ such that
 the two polynomials have the same largest zero  
$z_0$, then 
$$
z_0^3  - (a_1^2 + a_2^2 + a_3^2) z_0  = 2a_1a_2a_3 u
\quad \hbox{ and } \quad 
z_0^3  - (b_1^2 + b_2^2 + b_3^2) z_0  = 2b_1b_2b_3 u$$
so that 
$$z_0(z_0^2-(a_1^2+a_2^2+a_3^2))/(a_1a_2a_3)
= z_0(z_0^2-(b_1^2+b_2^2+b_3^2))/(b_1b_2b_3),$$
which can only happen for one positive value $z_0$ unless
$(\alpha_0,\alpha_1) = (\beta_0,\beta_1)$.
If $(\alpha_0,\alpha_1) \ne (\beta_0,\beta_1)$, then 
the difference of  the largest zeros of
$$
z_0^3  - (a_1^2 + a_2^2 + a_3^2) z_0  = 2a_1a_2a_3 u
\quad \hbox{ and } \quad 
z_0^3  - (b_1^2 + b_2^2 + b_3^2) z_0  = 2b_1b_2b_3 u$$
is always nonnegative  for $u \in [-1,1]$ if it is true 
for $u$ at the end points.

The remaining consequences can be readily verified. \qed

\medskip
Next, we consider the result for $n = 4$.

\begin{lemma} \label{lem} Suppose  
$R = S(r_1,r_2, r_3, r_4)$ with $r_1, r_2, r_3, r_4 > 0$. Then 
$$\det(zI - 2\Re(e^{i\theta}R)) =
z^4 - (r_1^2 + r_2^2 + r_3^2 + r_4^2)z^2 +
 r_1^2r_3^2 + r_2^2r_4^2 - 2r_1r_2r_3r_4\cos(4\theta).$$
\end{lemma}

\begin{theorem} \label{n=4} Let $A = S(a_1,\dots,a_4)$ and $B = S(b_1, \dots, b_4)$, 
where $\ba  = (a_1, a_2, a_3, a_4)$ and $\bb = (b_1, b_2, b_3, b_4)$ are nonnegative 
vectors. Suppose 
$$\alpha_1 = \sum_{j=1}^4 a_j^2, \ 
\alpha_2 = a_1^2a_3^2 + a_2^2a_4^2, \  \alpha_3 = \prod_{j=1}^4 a_j, \
\beta_1 = \sum_{j=1}^4 b_j^2, \
\beta_2 = b_1^2b_3^2 + b_2^2b_4^2, \  \beta_3 = \prod_{j=1}^4 b_j.$$
\begin{itemize}
\item[{\rm (a)}] We have $W(B) \subseteq W(A)$ if and only if
$$\alpha_1 + \sqrt{\alpha_1^2 -4(\alpha_2 \pm  2\alpha_3 )} 
\ge \beta_1 + \sqrt{\beta_1^2 -4(\beta_2  \pm 2 \beta_3)}.$$
\item[{\rm (b)}]
Suppose $(\alpha_1, \alpha_3) = (\beta_1, \beta_3)$. 
Then $W(B) \subseteq W(A)$ 
if and only if $a_1^2a_3^2+a_2^2a_4^2 \le b_1^2 b_3^2 + b_2^2 b_4^2$. 
\item[{\rm (c)}] Suppose
$a_1 \ge a_2 \geq a_3 \geq a_4 >0$ or $a_4 \geq a_3 \geq a_2 \geq a_1 >0$.  
Then 
$$W(S(a_1,a_3,a_2,a_4))  = W(S(a_1, a_4, a_2, a_3))$$
$$\subseteq 
 W(S(a_1,a_2,a_3,a_4)) = W(S(a_1, a_4, a_3, a_2))$$
$$ \subseteq 
W(S(a_1,a_2,a_4,a_3)) = W(S(a_1,a_3,a_4,a_2)).$$
\end{itemize}
\end{theorem}

{\it Proof.} (a) Note that $W(B) \subseteq W(A)$ if and only if  largest zero of 
$\det(\lambda I - 2 \Re(e^{i\theta}A))$ is larger than or equal to that 
of $\det(\lambda I - 2 \Re(e^{i\theta}B))$, i.e., 
 $$\alpha_1 + \sqrt{\alpha_1^2 -4(\alpha_2 - 2\alpha_3 u)} 
\ge \beta_1 + \sqrt{\beta_1^2 -4(\beta_2 -2 \beta_3 u)}$$
for any $u \in [-1,1]$. By elementary calculus, 
we see that the above is true for $u \in [-1,1]$ provided that  the inequality holds
for $u \in \{-1,1\}$.

(b) This follows readily from (a).

(c) 
By Lemma \ref{lem}, we get the set equalities.
To prove the set inclusions, we compare  $\lambda_1(\Re(e^{i\theta}A_\sigma))$
for different permutation $\sigma$. The result follows from  (b)
 and the following inequalities: 
\begin{eqnarray*}  
& & (a_1^2 a_3^2+a_2^2 a_4^2) - (a_1^2 a_4^2 +a_2^2 a_3^2) = (a_1^2-a_2^2)(a_3^2-a_4^2)
   \geq  0, \\ 
   & & (a_1^2 a_2^2+a_3^2 a_4^2) - (a_1^2 a_3^2 +a_2^2 a_4^2) = (a_1^2-a_4^2)(a_2^2-a_3^2)
   \geq 0. \
\end{eqnarray*} 
\vskip-.3in  \qed

\section{Result for $n = 5$}

In this section, we focus on the case that $n = 5$.

\begin{lemma} \label{p5} Let $B = S(b_1, \dots, b_5)$. Then 
$$\det(xI - (e^{i\theta}B+e^{-i\theta}B^t)) = x^5 - \beta_2 x^3 + \beta_1 x -
\beta_0 \cos(5 \theta)$$
with 
$$\beta_2 = b_1^2 + b_2^2 + b_3^2 + b_4^2 + b_5^2, \quad
\beta_1 = b_1^2b_3^2 + b_1^2b_4^2 + b_2^2b_4^2 + b_2^2b_5^2 + b_3^2b_5^2, \quad 
 \beta_0 = 2b_1b_2b_3b_4b_5.$$
\end{lemma}

\begin{proposition} \label{criterion} Let $B = S(b_1, \dots, b_5), 
C = S(c_1, \dots, c_5)\in M_5$, where 
$\bb = (b_1, \dots, b_5)$ and $\bc = (c_1, \dots, c_5)$ 
are nonnegative vectors such that
$\sum_{j=1}^5 b_j^2 = \sum_{j=1}^5 c_j^2$ and $b_1 \cdots b_5 = c_1\cdots c_5$.
Then $W(C) \subseteq W(B)$ if and only if
$$
b_1^2b_3^2 + b_1^2b_4^2 + b_2^2b_4^2 + b_2^2b_5^2 + b_3^2b_5^2
\le c_1^2c_3^2 + c_1^2c_4^2 + c_2^2c_4^2 + c_2^2c_5^2 + c_3^2c_5^2.$$
Furthermore, if $b_1^4 +b_2^4 +b_3^4 +b_4^4 +b_5^4=c_1^4 +c_2^4 +c_3^4 +c_4^4 +c_5^4$, 
the criterion inequality is equivalent to 
\begin{equation}\label{c1c2}
c_1^2 c_2^2 +c_2^2 c_3^2 +c_3^2 c_4^2 +c_4^2 c_5^2 +c_5^2 c_1^2 \leq 
b_1^2 b_2^2 +b_2^2 b_3^2 +b_3^2 b_4^2 +b_4^2 b_5^2 +b_5^2 b_1^2.
\end{equation}
\end{proposition}

\it Proof. \rm The proof follows from Lemmas \ref{poly}, \ref{p5}, and the fact that
 $$(b_1^2 +b_2^2 +b_3^2 +b_4^2 +b_5^2)^2 -(b_1^4 +b_2^4 +b_3^4 +b_4^4 +b_5^4)$$
 $$=2(b_1^2 b_2^2 +b_2^2 b_3^2 +b_3^2 b_4^2 +b_4^2 b_5^2 +b_5^2 b_1^2) +2(b_1^2 b_3^2 +b_1^2 b_4^2 +b_2^2 b_4^2 +b_2^2 b_5^2 +b_3^2 b_5^2).$$
\vskip -.3in \qed

\medskip

 \begin{theorem} \label{main-5} Suppose $n = 5$ and $a_1 \geq a_2 \geq a_3 \geq a_4 \geq a_5 > 0$. 
 Then for $\sigma\in S_5$,
$$W(S(a_5, a_1, a_4, a_3, a_2)) \subseteq
W(S(a_{\sigma(1)}, a_{\sigma(2)}, \ldots, a_{\sigma(5)})) 
\subseteq W(S(a_1,a_3,a_5,a_4,a_2)).$$
\end{theorem}  

\it Proof. \rm 
By  Proposition \ref{criterion}, 
we only need to compare the cyclic sums:
$$a_{\sigma(1)}^2 a_{\sigma(2)}^2 +a_{\sigma(2)}^2 a_{\sigma(3)}^2 +a_{\sigma(3)}^2 a_{\sigma(4)}^2 +a_{\sigma(4)}^2 a_{\sigma(5)}^2 +a_{\sigma(5)}^2 a_{\sigma(1)}^2$$
in (\ref{c1c2}) for all permutations $\sigma \in S_5$.

It is known (\cite{CN1}) that the numerical range $W(S(a_1, a_2, a_3, a_4, a_5))$ is invariant under a cyclic permutation (i.e.  $W(S(a_1, a_2, a_3, a_4, a_5))
=W(S(a_2, a_3, a_4, a_5, a_1))$), and a reversal permutation \\ (i.e. $W(S(a_1, a_2, a_3, a_4, a_5))=W(S(a_5, a_4, a_3, a_2, a_1))$).  
Hence, comparing the cyclic sums for all permutations $\sigma \in S_5$, it remains to compare only the following 12 permutations:
\medskip

 (1) $S(1)=  S(a_3, a_5, a_4, a_2, a_1)$,
\hskip .6in 
 (2) $S(2)= S(a_3, a_5, a_4, a_1, a_2)$,

 (3) $S(3)= S(a_2,a_5, a_4,a_1,a_3)$, 
\hskip .6in
 (4) $S(4)= S(a_2, a_5, a_4, a_3,a_1)$,

(5) $S(5)=S(a_2,a_5, a_3, a_4, a_1)$, 
\hskip .6in (6) $S(6)=S(a_2, a_5, a_3, a_1, a_4)$, 

(7) $S(7)=S(a_1, a_5, a_4, a_3, a_2)$, 
\hskip .6in (8) $S(8)=S(a_1, a_5, a_4, a_2, a_3)$, 

(9) $S(9)=S(a_1, a_5, a_3, a_4, a_2)$, 
\hskip .6in (10) $S(10)=S(a_1, a_5, a_3, a_2, a_4)$, 

(11) $S(11)=S(a_1, a_5, a_2, a_4, a_3)$, 
\hskip .45in (12) $S(12)=S(a_1, a_5, a_2, a_3, a_4)$. 

We  introduce a partial order
 $S(b_1,b_2,b_3, b_4, b_5)  \geq S(c_1, c_2, c_3, c_4, c_5)$ by the inequality 
 $$b_1^2 b_2^2 +b_2^2 b_3^2 +b_3^2 b_4^2 +b_4^2 b_5^2 +b_5^2 b_1^2 \geq c_1^2 c_2^2 +c_2^2 c_3^2 +c_3^2 c_4^2 +c_4^2 c_5^2 +
 c_5^2 c_1^2.$$ 
For ease of notation,  the cyclic sum of $S(j)$ is also denoted by by $S(j)$.

Let
$a_1^2 =r_1, a_2^2 =r_1 +r_2, a_3^2 =r_1 +r_2 +r_3, a_4^2=r_1 +r_2 +r_3 +r_4, a_5^2 =r_1 +r_2 +r_3 +r_4 +r_5$. 
Direct computations show that
 \begin{eqnarray*}
&& S(1) -S(2) =r_2 r_4 \geq 0, \ S(2) -S(3) =r_3 (r_2 +r_3 +r_4 +r_5) \geq 0,\\
&& S(3) -S(8) =r_2 r_5 \geq 0, \  S(8) -S(10)=r_4 (r_2 +r_3 +r_4 +r_5) \geq 0,\\
&& S(10)-S(12)=r_3 r_5 \geq 0,\ S(5) -S(6) =r_3 (r_2 +r_3 +r_4) \geq 0,\\ 
&& S(1)- S(4) =r_3 r_5 \geq 0, \ S(4) -S(7) =r_2 (r_4 +r_5) \geq 0,\\
&& S(7) -S(8) =r_3 (r_2 +r_3 +r_4) \geq 0,\ S(4) - S(5)=r_4 (r_2 +r_3 +r_4 +r_5) \geq 0,\\
&& S(5) -S(9)=r_2 r_5 \geq 0, \ S(9)-S(11)=r_3 (r_2 +r_3 +r_4 +r_5) \geq 0,\\
&& S(11)-S(12)=r_2 r_4 \geq 0, \  S(6)-S(10)=r_2 (r_4 +r_5) \geq 0.
\end{eqnarray*}
These relations imply that $S(12)\le S(j)\le S(1)$ for all $j=1,2,\ldots 12$, and thus the numerical range inclusions  follow. \qed

\begin{remark} \rm
We may apply a result of Chang and Wang \cite[Theorem 2.1]{CW} to give an alternative proof for the right hand inclusion of Theorem \ref{main-5}. Let $A_\sigma = S(a_{\sigma(1)}, a_{\sigma(2)},\dots, a_{\sigma(5)})$ and  $f_\theta(x)=\det(xI - (e^{i\theta}A_\sigma+e^{-i\theta}A_\sigma^t))$.
Then by Lemma \ref{p5}
$$f_\theta(x)=f_0(x)-2a_{\sigma(1)}a_{\sigma(2)}\cdots a_{\sigma(5)}\cos(5\theta).$$
Similarly, let  $A_\tau=S(a_1,a_3,a_5,a_4,a_2)$ and $g_\theta(x)=\det(xI - (e^{i\theta}A_\tau+e^{-i\theta}A_\tau))$. We have
$$g_\theta(x)=g_0(x)-2a_{\tau(1)}a_{\tau(2)}\cdots a_{\tau(5)}\cos(5\theta).$$
Then, by  \cite[Theorem 2.1]{CW},
$$f_\theta(x)-g_\theta(x)=f_0(x)-g_0(x)\ge 0,$$
and thus the largest zero of $g_\theta(x)$ is larger than that of $f_\theta(x)$.
Therefore, $\lambda_1(\Re(e^{i\theta}A_\tau))\ge \lambda_1(\Re(e^{i\theta}A_\sigma))$ for any $\sigma \in S_5$, and so $W(S(a_{\sigma(1)}, a_{\sigma(2)}, \ldots, a_{\sigma(5)})) 
\subseteq W(S(a_1,a_3,a_5,a_4,a_2))$.
\end{remark}

\section{Results for $n=6$} 

Let $a_6 \ge \cdots \ge a_1 \ge 0$. 
In the following, we consider
$A _\sigma = S(a_{\sigma(1)},a_{\sigma(2)}, \ldots, a_{\sigma(6)})$ 
for $\sigma \in S_6$, and determine a permutation such that
$W(A_\sigma)$ is maximum in terms of set inclusion. The following is our main theorem.

 \begin{theorem} \label{upper=bound45678} Suppose $(a_1,a_2, a_3,a_4,a_5,a_6)$  has positive entries arranged in ascending or descending order.
  For any $\sigma \in S_6$, the following inclusion relation holds:
  $$W(S(a_{\sigma(1)}, a_{\sigma(2)}, a_{\sigma(3)}, a_{\sigma(4)}, a_{\sigma(5)}, a_{\sigma(6)})) \subseteq W(S(a_1,a_3,a_5,a_6,a_4,a_2)).$$
\end{theorem}

\medskip 

Theorem \ref{upper=bound45678} will be proved by a series of lemmas and propositions. 
Firstly, we have the following.

\medskip  
\begin{lemma}\label{det}
 Let $A=S(a_1,a_2,a_3,a_4,a_5,a_6)$. Then
\begin{eqnarray}\label{detpoly}
&&\det(xI - e^{i\theta}A - e^{-i\theta}A^*)\nonumber\\
&=&x^6 - (a_1^2 + a_2^2 + a_3^2 + a_4^2 + a_5^2 + a_6^2)x^4\nonumber\\
&&+(a_1^2 a_3^2 + a_1^2 a_4^2 + a_1^2 a_5^2 + a_2^2 a_4^2  + a_2^2 a_5^2 + a_2^2 a_6^2 + a_3^2 a_5^2 + a_3^2 a_6^2 + a_4^2 a_6^2)x^2 \\ 
&&-a_1^2 a_3^2 a_5^2 - a_2^2 a_4^2 a_6^2 - 2a_1 a_2 a_3 a_4 a_5 a_6\cos(6\theta).\nonumber
\end{eqnarray}
\end{lemma}

\medskip

The coefficient $S_2(a_1,\ldots, a_6)=a_1^2 a_3^2 +\ldots +a_4^2 a_6^2$
in (\ref{detpoly}) and the cyclic sum $\sum_{i=1}^6 a_i^2 a_{i+1}^2$ 
of $(a_1, \ldots,a_6)$ are related by the equation:
\begin{equation}\label{equi} 2\Big(\sum_{1 \leq i <j \leq 6, 2 \leq j-i \leq 5} a_i^2 a_j^2 +\sum_{i=1}^6 a_i^2 a_{i+1}^2 \Big)=\left(\sum_{i=1}^6 a_i^2\right)^2 -\sum_{i=1}^6 a_i^4,\end{equation} 
where $a_7=a_1$.

By replacing $x^2$ by $x$, the determinantal polynomial (\ref{detpoly}) becomes 
\begin{equation}\label{cubic}
\begin{aligned}
F_A(x: \theta)& =x^3 -(a_1^2 +a_2^2 +a_3^2 +a_4^2 +a_5^2 +a_6^2) x^2 +(a_1^2 a_3^2 +a_1^2 a_4^2+\cdots +a_4^2 a_6^2) x \\
  &\quad -a_1^2 a_3^2 a_5^2 -a_2^2 a_4^2 a_6^2 -2 a_1 a_2 a_3 a_4 a_5 a_6 \cos(6\theta).
\end{aligned}
\end{equation}
Let $\sigma\in S_6$ and $A=S(a_{\sigma(1)}, a_{\sigma(2)},\ldots,a_{\sigma(6)})$. In this case, the cubic polynomial $F_A(x: \theta)$ is also denoted by $F_{\sigma}(x: \theta)$. Indeed,  
\begin{equation*}
\begin{aligned}
F_{\sigma}(x: \theta)& =x^3 -(a_1^2 +a_2^2 +a_3^2 +a_4^2 +a_5^2 +a_6^2) x^2 +(a_{\sigma(1)}^2 a_{\sigma(3)}^2 +a_{\sigma(1)}^2 a_{\sigma(4)}^2+\cdots +a_{\sigma(4)}^2 a_{\sigma(6)}^2) x \\
  &\quad -a_{\sigma(1)}^2 a_{\sigma(3)}^2 a_{\sigma(5)}^2 -a_{\sigma(2)}^2 a_{\sigma(4)}^2 a_{\sigma(6)}^2 -2 a_1 a_2 a_3 a_4 a_5 a_6 \cos(6\theta).
\end{aligned}
\end{equation*}
The maximal eigenvalue problem of $\Re(e^{i\theta}A)$ is reduced to  the maximal zero of the cubic polynomial $F_{\sigma}(x: \theta)$ over $0\le \theta< 2\pi$.
For two permutations $\sigma, \eta \in S_6$, we define and denote  $S(a_{\eta(1)}, \ldots, a_{\eta(6)}) \leq S(a_{\sigma(1)}, \ldots, a_{\sigma(6)})$ if  their cyclic sums satisfy
  $$\sum_{i=1}^6 a_{\eta(i)}^2 a_{\eta(i+1)}^2 \leq \sum_{i=1}^6 a_{\sigma(i)}^2 a_{\sigma(i+1)}^2.$$
By (\ref{equi}), 
$S(a_{\eta(1)}, \ldots, a_{\eta(6)}) \leq S(a_{\sigma(1)}, \ldots, a_{\sigma(6)})$ 
if and only if 
$$\sum_{1 \leq i <j \leq 6, 2 \leq j-i \leq 5} a_{\sigma(i)}^2 a_{\sigma(j)}^2
\le \sum_{1 \leq i <j \leq 6, 2 \leq j-i \leq 5} a_{\eta(i)}^2 a_{\eta(j)}^2.$$
The cyclic sum $\sum_{i=1}^6 a_{\sigma(i)}^2 a_{\sigma(i+1)}^2$ plays an important role in determining the maximal zero of $F_{\sigma}(x,\theta))$. 

\medskip
It is clear that the polynomial $F_{\sigma}(x:\theta)$ is invariant under the cyclic permutation \\ $(\sigma(2),\sigma(3),\sigma(4),\sigma(5),\sigma(6),\sigma(1))$
and the reversal permutation $(\sigma(6),\sigma(5),\sigma(4),\sigma(3),\sigma(2),\sigma(1))$   of $\sigma$. The $6!=720$ permutations $\sigma \in S_6$ have $60$ representatives under this invariance. 
The $60$ representative cyclic shift matrices are classified into $10$ families (I) to (X)  according to the  equivalence relation $\sigma \cong \eta$ defined by 
  $$\{\sigma(1), \sigma(3), \sigma(5)\}=\{\eta(1), \eta(3), \eta(5)\}, \ {\rm equivalently,}\ \{\sigma(2), \sigma(4), \sigma(6)\}=\{\eta(2), \eta(4), \eta(6)\}.$$
We explicitly list $60$ representatives $S(k)$, $k=1,\ldots,60$, in the respective 10 families.

\begin{itemize}
\item[(I)] 
 $S(21)=S(a_4,a_2, a_1, a_3,a_5, a_6)$,
 $S(28)=S(a_1, a_3, a_5, a_2, a_4, a_6)$,
 $S(2)=S(a_1, a_2, a_4, a_3, a_5, a_6)$, $S(4)=S(a_1, a_3, a_4, a_2, a_5, a_6)$,
$S(23)=S(a_4,a_3, a_1,a_2, a_5, a_6), S(26)=S(a_1, a_2, a_5, a_3, a_4, a_6)$.
\item[(II)] 
 $S(45)=S(a_2, a_4, a_1, a_5, a_3, a_6)$, $S(53)=S(a_1, a_5, a_3, a_4, a_2, a_6)$, $S(47)=S(a_2, a_5, a_1, a_3, a_6)$, 
$S(51)=S(a_1, a_4, a_3, a_5, a_2, a_6)$, $S(57)=S(a_1, a_4, a_2, a_5, a_3, a_6)$,
  $S(59)=S(a_1, a_5, a_2, a_4, a_3, a_6)$. 
\item[(III)]
 $S(15)=S(a_3,a_2, a_1, a_4, a_5, a_6), S(58)=S(a_1, a_4, a_5, a_2, a_3, a_6)$, $S(1)=S(a_1,a_2, a_3, a_4, a_5, a_6)$,\\ $S(6)=S(a_1, a_4, a_3, a_2, a_5, a_6)$, $S(17)=S(a_3, a_4, a_1, a_2, a_5, a_6)$,
 $S(56)=S(a_1, a_2, a_5, a_4, a_3, a_6)$.
\item[(IV)] 
$S(9)=S(a_2, a_3, a_1, a_4, a_5, a_6), S(52)=S(a_1, a_4, a_5, a_3, a_2, a_6)$,
$S(3)=S(a_1,a_3, a_2, a_4, a_5, a_6)$,\\ $S(5)=S(a_1, a_4, a_2, a_3, a_5, a_6)$, 
$S(11)=S(a_2,a_4, a_1,a_3, a_5, a_6)$, 
$S(50)=S(a_1, a_3, a_5, a_4, a_2, a_6)$.
\item[(V)] 
 $S(13)=S(a_3, a_1, a_2, a_4, a_5, a_6), S(46)=S(a_2, a_4, a_5, a_1, a_3, a_6)$,
 $S(7)=S(a_2, a_1, a_3, a_4, a_5, a_6)$,\\ 
$S(12)=S(a_2, a_4, a_3, a_1, a_5, a_6)$, $S(18)=S(a_3, a_4, a_2, a_1, a_5, a_6)$, 
 $S(44)=S(a_2, a_1, a_5, a_4, a_3, a_6)$.
\item[(VI)] 
 $S(19)=S(a_4, a_1, a_2, a_3, a_5, a_6),S(34)=S(a_2, a_3, a_5, a_1, a_4, a_6)$, $S(8)=S(a_2, a_1, a_4, a_3, a_5, a_6)$,\\
 $S(10)=S(a_2, a_3, a_4, a_1, a_5,a_6)$, $S(24)=S(a_4, a_3, a_2, a_1, a_5, a_6)$, 
 $S(32)=S(a_2, a_1, a_5, a_3, a_4)$.
\item[(VII)] 
$S(20)=S(a_4, a_1, a_3, a_2, a_5, a_6), S(40)=S(a_3, a_2, a_5, a_1, a_4, a_6)$,
 $S(14)=S(a_3, a_1, a_4, a_2, a_5, a_6)$,\\ $S(16)=S(a_3, a_2, a_4, a_1, a_5, a_6)$, $S(22)=S(a_4, a_2, a_3, a_1, a_5, a_6)$, 
 $S(38)=S(a_3, a_1, a_5, a_2, a_4, a_6)$. 
\item[(VIII)]
 $S(39)=S(a_3, a_2, a_1, a_5, a_4, a_6), S(60)=S(a_1, a_5, a_4, a_2, a_3, a_6)$,
  $S(25)=S(a_1, a_2, a_3, a_5, a_4, a_6)$,\\ $S(30)=S(a_1, a_5, a_3, a_2, a_4, a_6), S(41)=S(a_3, a_5, a_1, a_2, a_4, a_6)$,  
  $S(55)=S(a_1, a_2, a_4, a_5, a_3, a_6)$. 
\item[(IX)] 
  $S(33)=S(a_2, a_3, a_1,a_5, a_4, a_6), S(54)=S(a_1, a_5, a_4, a_3, a_2, a_6)$,
  $S(27)=S(a_1, a_3, a_2, a_5, a_4, a_6)$,\\ $S(29)=S(a_1, a_5, a_2, a_3, a_4, a_6)$,  $S(35)=S(a_2, a_5, a_1, a_3, a_4, a_6)$,
 $S(49)=S(a_1, a_3, a_4, a_5, a_2, a_6)$.
\item[(X)] 
 $S(37)=S(a_3, a_1, a_2, a_5, a_4, a_6),S(48)=S(a_2, a_5, a_4, a_1, a_3, a_6)$,
 $S(31)=S(a_2, a_1, a_3, a_5, a_4, a_6)$,\\ $S(36)=S(a_2, a_5, a_3, a_1, a_4, a_6)$, $S(42)=S(a_3, a_5, a_2, a_1, a_4, a_6)$,
 $ S(43)=S(a_3, a_5, a_4, a_1, a_2, a_6)$. 
\end{itemize}

 \medskip 

 For two permutations $\sigma, \eta \in S_6$ which are neither  cyclic equivalent nor reversal equivalent,  
the inclusion $W(S(a_{\eta(1)},\ldots, a_{\eta(6)})) \subseteq W(S(a_{\sigma(1)}, \ldots, a_{\sigma(6)})$ holds if and only if the greatest zero of $F_{\eta}(x:\theta)$ is less than 
or equal to the greatest zero of $F_{\sigma}(x: \theta)$ for any angle $\theta \in [0, 2\pi]$. 
The comparison falls into two types: 

\medskip
\quad $(i)$ $\sigma$ and $\eta$ are in the same family;
\qquad 
$(ii)$  $\sigma$ and $\eta$ are in distinct families.  

\medskip
For the comparison of type $(i)$, we need the following result.

\begin{proposition}\label{prop3} 
Let $\bc=(c_1, c_2, c_3)$ and $\bd=(d_1, d_2, d_3)$ be positive vectors with descending entries such that 
  $c_1 +c_2 +c_3 =d_1+ d_2+ d_3$ and $c_1 c_2 c_3 =d_1 d_2 d_3$.
If $c_1 c_2 +c_1 c_3 +c_2 c_3 < d_1 d_2 +d_1 d_3 +d_2 d_3$, then $d_1 < c_1$. 
\end{proposition}

 {\it Proof.} Let $ \alpha=c_1+c_2 +c_3, \ \beta= c_1 c_2 c_3, \ a = d_1 d_2 +d_1 d_3 +d_2 d_3$ and $\tilde{a}=c_1 c_2 +c_1 c_3 +c_2 c_3$. Then the two cubic polynomials
  $$P(t)=t^3 -\alpha t^2 +a t -\beta,  \quad \tilde{P}(t)=t^3 -\alpha t^2 +\tilde{a} t -\beta$$
 satisfy $P(t)-\tilde{P}(t)=(a -\tilde{a}) t >0$ for all $t >0$. Since $c_1 >0$, it implies $P(t)-\tilde{P}(t) >0$ for all $t \geq c_1$. The assertion follows from Lemma 2.4 (b). \qed

Applying Proposition \ref{prop3}, we easily obtain the following criterion for type $(i)$.

 \begin{proposition} \label{n=6cri}
 Let $\bb=(b_1, b_2, b_3, b_4, b_5, b_6) \in {\mathbb R}^6$ and $\bc=(c_1, c_2, c_3, c_4, c_5, c_6) \in {\mathbb R}^6$ be positive vectors 
 such that $\bb \bb^t=\bc \bc^t$ and
  $$ b_1 b_2 b_3 b_4 b_5 b_6=c_1 c_2 c_3 c_4 c_5 c_6,\quad 
  b_1^2 b_3^2 b_5^2 +b_2^2 b_4^2 b_6^2 =c_1^2 c_3^2 c_5^2 +c_2^2 c_4^2 c_6^2,$$
  $$  b_1^4 +b_2^4 +b_3^4 +b_4^4 +b_5^4 +b_6^4=c_1^4 +c_2^4 +c_3^4 +c_4^4 +c_5^4 +c_6^4. $$
 If the cyclic sums satisfy the inequality
\begin{equation}\label{bccr} c_1^2 c_2^2 +c_2^2 c_3^2 +c_3^2 c_4^2 +c_4^2 c_5^2 +c_5^2 c_6^2 +c_6^2 c_1^2 \leq b_1^2 b_2^2 +b_2^2 b_3^2 +b_3^2 b_4^2 +b_4^2 b_5^2 +b_5^2 b_6^2 +b_6^2 b_1^2
\end{equation} 
then  $W(S(c_1, c_2, c_3, c_4, c_5, c_6)) \subseteq W(S(b_1, b_2, b_3, b_4, b_5, b_6)).$
\end{proposition}

\it Proof. \rm Let $B=S(b_1, \dots, b_6)$ and $C=S(c_1, \dots, c_6)$. 
Applying Proposition \ref{prop3} to the two cubic polynomial $F_B(x:\theta)$ and $F_C(x:\theta)$ of the form (\ref{cubic}), if  
\begin{equation}\label{bc}
b_1^2b_3^2 + b_1^2b_4^2 + \cdots + b_3^2b_5^2 
+ b_3^2b_6^2 + b_4^2b_6^2\\
\le c_1^2c_3^2 + c_1^2c_4^2 + c_1^2c_5^2 + \cdots + c_3^2c_5^2 
+ c_3^2c_6^2 + c_4^2b_c^2
\end{equation}
then $\lambda_1(\Re(e^{i\theta}C))\le \lambda_1(\Re(e^{i\theta}B))$, 
and thus $W(S(c_1,\dots,c_6))\subset W(S(b_1,\dots,b_6))$.
Now, according to the identity (\ref{equi}), the condition (\ref{bc}) 
is equivalent to inequality (\ref{bccr}). 
\qed

\medskip

Again, we define a partial order
 $S(b_1,b_2,b_3, b_4, b_5,b_6)  \le S(c_1, c_2, c_3, c_4, c_5,c_6)$ according to their cyclic sums
 $$b_1^2 b_2^2 +b_2^2 b_3^2 +b_3^2 b_4^2 +b_4^2 b_5^2 +b_5^2 b_6^2+b_6^2b_1^2 \le c_1^2 c_2^2 +c_2^2 c_3^2 +c_3^2 c_4^2 +c_4^2 c_5^2 +
 c_5^2 c_6^2+c_6^2c_1^2.$$ 
For notational convenience, if no confusion is caused, the  
cyclic sum of the cyclic shift matrix $S(j)$ in the ten different 
classes (I) -- (X) will also be denoted by  $S(j)$.

For the cyclic shift matrix $S(a_1, a_2, a_3, a_4, a_5, a_6)$ with positive weights, we 
let
  $$a_1^2=r_1,\ a_2^2=r_1 +r_2,\ a_3^2=r_1+r_2 +r_3,\ a_4^2=r_1+r_2+r_3+r_4,$$
 $$ a_5^2=r_1+r_2+r_3+r_4+r_5,\ a_6^2=r_1+r_2+r_3+r_4+r_5+r_6.$$
 If $0 < a_1 < a_2 <a_3 < a_4 < a_5 < a_6$ then $r_1,\ldots r_6$ are positive.

\medskip 
 \begin{lemma} \label{top21}  Let $0 < a_1 < a_2 <a_3 < a_4 < a_5 < a_6$. With respect to the cyclic sum the following inequalities hold. 
 
\begin{description}
 \item[{\rm (I)}] $S(j) \leq S(21)$ for $j=2,4,23,26,28$,
 \item[{\rm (II)}]  $S(j) \leq S(45)$ for $j=47,51,53,57,59$, 
 \item[{\rm (III)}]  $S(j) \leq S(15)$ for $j=1,6,17,56,58$,
 \item[{\rm (IV)}]  $S(j) \leq S(9)$  for $j=3,5,11,50,52$, 
 \item[{\rm (V)}]  $S(j) \leq S(13)$  for $j=7,12,18,44,46$,
 \item[{\rm (VI)}]  $S(j) \leq S(19)$ for $j=8,10,24,32,34$,
 \item[{\rm (VII)}] $S(j) \leq S(20)$ for $j=14,16,22,38,40$,
  \item[{\rm (VIII)}] $S(j) \leq S(39)$  for $j=25,30,41,55,60$,
 \item[{\rm (IX)}]  $S(j) \leq S(33)$  for $j=27,29,35,49,54$,
  \item[{\rm (X)}] $S(j) \leq S(37)$ for $j=31, 36,42,43,48$.
 \end{description}
Moreover the inequalities also hold: 
\begin{description}
  \item[\rm (a)]   $S(45) \leq S(33) \leq S(37) \leq S(13) \leq S(21)$, 
  \item[\rm (b)] $S(33) \leq S(20) \leq S(19) \leq S(21)$,
  \item[\rm (c)] $S(33) \leq S(39) \leq S(19) \leq S(21)$,
  \item[\rm (d)] $S(9) \leq S(15) \leq S(13) \leq S(21)$.
\end{description} 
\end{lemma} 
   
\medskip 
\it Proof. \rm We compute the differences of cyclic sums of $S(j)$'s in each family. 
\begin{itemize}
\item[{\rm (I)}] The differences $S(21)-S(2), S(21)-S(4), S(21)-S(23),S(21)-S(26),S(21)-S(28)$ are respectively computed by
  $(r_2 +r_3 +r_4)(r_4 +r_5 +r_6), r_3 (r_2 +r_3 +r_4 +r_5) +(r_2 +r_3 +r_4)(r_4 +r_5 +r_6)$, $r_3 r_5 ,r_5 (r_3 +r_4 +r_5 +r_6)+
(r_2 +r_3 +r_4)(r_4 +r_5 +r_6),(r_2 +r_3 +r_4 +r_5)(r_3 +r_4 +r_5 +r_6)$.
\item[{\rm (II)}] 
The differences $S(45)-S(47), S(45)-S(51), S(45)-S(53), S(45)-S(57), S(45)-S(59)$ are respectively expressed as $r_3 r_5, \,r_3 r_5 +(r_2 +r_3) r_6,
(r_2 +r_3)(r_5 +r_6)$, $r_2 r_6,(r_2 +r_3) r_5 +r_2 r_6$.
\item[{\rm (III)}] The differences $S(15)-S(1), S(15)-S(6), S(15)-S(17)$, $S(15)-S(56),S(15)-S(58)$ are respectively given by   
 $(r_2 +r_3)(r_5 +r_6),(r_2+r_3)(r_5+r_6) +(r_3 +r_4)(r_2 +r_3+r_4+r_5)$, $(r_3+r_4)(r_4+r_5), (r_3 +r_4)(r_4 +r_5) +(r_2+r_3+r_4+r_5)(r_5+r_6),
  (r_2 +r_3 +r_4 +r_5)(r_3 +r_4+r_5 +r_6)$.
\item[{\rm (IV)}] 
The differences $S(9)-S(3), S(9)-S(5), S(9)-S(11)$, $S(9)-S(50),S(9)-S(52)$ are respectively computed by  $r_2 (r_5 +r_6), r_2 (r_5 +r_6) +
r_4 (r_2 +r_3+r_4 +r_5), r_4 (r_3 +r_4 +r_5)$, $r_4 (r_3+r_4 +r_5) +(r_5 +r_6)(r_2+r_3+r_4+r_5),(r_4+r_5+r_6)(r_2+r_3+r_4+r_5)$.
\item[{\rm (V)}] 
The differences $S(13)-S(7), S(13)-S(12), S(13)-S(18)$, $S(13)-S(44), S(13)-S(46)$ are respectively given by 
 $r_3 (r_5 +r_6), r_3 (r_5 +r_6) +(r_2 +r_3 +r_4)(r_3+r_4 +r_5)$, $(r_4 +r_5)(r_2+r_3 +r_4),(r_4 +r_5)(r_2 +r_3 +r_4) +(r_5+r_6)(r_3+r_4+r_5),(r_3+r_4+r_5)(r_2+r_3+r_4+r_5+r_6)$. 
\item[{\rm (VI)}] 
The differences $S(19)-S(8), S(19)-S(10), S(19)-S(24)$, $S(19)-S(32), S(19)-S(34)$ are respectively given by $(r_3 +r_4)(r_4+r_5+r_6), (r_2 +r_3)(r_3+r_4+r_5) +(r_3+r_4)(r_4+r_5+r_6), r_5 (r_2 +r_3)$, $r_5 (r_2 +r_3) +(r_3 +r_4 +r_5)(r_4 +r_5 +r_6), (r_3 +r_4 +r_5)(r_2+r_3+r_4+r_5 +r_6)$.
 \item[{\rm (VII)}] 
The differences $S(20)-S(14), S(20)-S(16), S(20)-S(22)$, $S(20)-S(38), S(20)-S(40)$ are respectively computed by
$r_4 (r_3+r_4 +r_5 +r_6), r_2 (r_4 +r_5) +r_4 (r_3+r_4+r_5 +r_6), r_2 r_5$, $r_4 (r_3 +r_4 +r_5 +r_6) +r_5 (r_2 +r_3 +r_4 +r_5 +r_6),
(r_4 +r_5)(r_2 +r_3 +r_4 +r_5 +r_6)$.
 \item[{\rm (VIII)}] 
The differences $S(39)-S(25),S(39)-S(30),S(39)-S(41)$, $S(39)-S(55), S(39)-S(60)$ are respectively computed by
 $r_6 (r_2 +r_3), r_6 (r_2 +r_3) +(r_2 +r_3 +r_4)(r_3+r_4+r_5), r_4 (r_3 +r_4 +r_5)$, $r_4 (r_3 +r_4 +r_5) +r_6 (r_2 +r_3 +r_4), (r_2 +r_3 +r_4)(r_3 +r_4 +r_5 +r_6)$.
\item[{\rm (IX)}]  
The differences $S(33)-S(27), S(33)-S(29), S(33)-S(35)$, $S(33)-S(49),S(33)-S(54)$ are respectively given by
$r_2 r_6, r_2 (r_4 +r_5 +r_6) +(r_3 +r_4)(r_4 +r_5),(r_3 +r_4)(r_4 +r_5)$. $(r_3 +r_4)(r_4 +r_5) +r_6 (r_2 +r_3 +r_4),(r_2 +r_3 +r_4)(r_4+r_5 +r_6)$.
\item[{\rm (X)}] 
The differences $S(37)-S(31), S(37)-S(36), S(37)-S(42)$, $S(37)-S(43), S(37)-S(48)$ are  respectively given by 
 $r_3 r_6, r_3 r_6 +(r_3+r_4)(r_2 +r_3 +r_4 +r_5), r_4 (r_2 +r_3 +r_4 +r_5)$, $r_3 r_6 +r_4(r_2+r_3+r_4+r_5+r_6), (r_3 +r_4)(r_2 +r_3 +r_4 +r_5 +r_6)$. 
\end{itemize}
All of the above differences of cyclic sums are positive, and thus the inequalities of 
(I)--(X) follow.
Furthermore, to prove (a), note that

\medskip
\quad  $S(33)-S(45)=r_4 (r_2 +2r_3+2r_4+2r_5+r_6)$, 

\medskip
\quad $S(37)-S(33)=r_2 (r_4+r_5)+r_3(r_2+r_3+r_4+r_5+r_6)$, 

\medskip
\quad $S(13)-S(37)=r_5 (r_3 +r_4 +r_5 +r_6)$, \quad
$S(21)-S(13)=r_4 r_6$.

\medskip\noindent
To prove (b), note that 

\medskip\quad
$ S(20)-S(33)=r_5 (r_2 +r_3 +r_4 +r_5) +r_6 (r_3 +r_4 +r_5)$, 

\medskip\quad
$ S(19)-S(20)=r_3 (r_2 +r_3 +r_4 +r_5)$, \quad $S(21)-S(19)=r_2 r_4$.

\medskip\noindent
To prove (c), note that

\medskip\quad
 $S(39) -S(33)=r_3 (r_2+r_3+r_4+r_5 +r_6)$, 
 
\medskip\quad 
 $S(19)-S(39)=r_4 r_6 +r_5 (r_2 +r_3 +r_4 +r_5 +r_6)$, \quad $S(21)-S(19)=r_2 r_4$.

\medskip\noindent
Finally, for (d), we have

\medskip
\quad  $S(15)-S(9)=r_3 (r_2 +r_3 +r_4+r_5+r_6)$, 

\medskip\quad $ S(13)-S(15)
=r_2 r_4$, \quad $S(21)-S(13)=r_4 r_6$.
\qed

\begin{remark} 
From Lemma \ref{top21}, we conclude that $S(j) \leq S(21)$ for all $j=1,2, \ldots, 60$.
\end{remark}

\medskip

By Proposition \ref{n=6cri} and Lemma \ref {top21}, we have the following.

\begin{proposition} \label{reduction}
        Each family of the 10 families (I)-(X) of the cyclic shift matrices has the largest numerical range. 
   \begin{description}
 \item{\rm (I)} $W(S(j)) \subseteq W(S(21))$  for $j=2,4,23,26,28$,
 \item{\rm (II)} $W(S(j)) \subseteq W( S(45))$  for $j=47,51,53,57,59$, 
 \item{\rm (III)} $W(S(j)) \subseteq W( S(15)) $  for $j=1,6,17,56,58$,
 \item{\rm (IV)} $W(S(j)) \subseteq  W( S(9))$    for $j=3,5,11,50,52$, 
\item{\rm (V)} $W(S(j)) \subseteq W(S(13))$ for $j=7,12,18,44,46$,
\item{\rm (VI)} $W(S(j)) \subseteq W(S(19))$ for $j=8,10,24,32,34$,
\item{\rm (VII)} $W(S(j)) \subseteq W( S(20))$  for $j=14,16,22,38,40$,
\item{\rm  (VIII)} $W(S(j)) \subseteq W S(39))$  for $j=25,30,41,55,60$,
\item{\rm  (IX)} $W(S(j)) \subseteq W( S(33))$  for $j=27,29,35,49,54$,
\item{\rm (X)} $W(S(j)) \subseteq W(S(37))$ and  for $j=31, 36,42,43,48$.
 \end{description}
 \end{proposition}

\medskip

Next, we deal with cyclic shift matrices with distinct cyclic sums belonging to type $(ii)$.

\medskip

\begin{proposition} \label{cubic3}
Let $\bc=(c_1, c_2, c_3)$ and $\bd=(d_1, d_2, d_3)$ be positive vectors with descending coordinates such that 
  $\alpha=c_1 +c_2 +c_3 =d_1+ d_2+ d_3$ and $c_1 c_2 c_3 \ne d_1 d_2 d_3$.
If $c_1 c_2 +c_1 c_3 +c_2 c_3 < d_1 d_2 +d_1 d_3 +d_2 d_3$ and the value $x_0$ determined by 
  $$\Big((d_1 d_2 +d_1 d_3 +d_2 d_3) -(c_1 c_2 +c_1 c_3 +c_2 c_3)\Big) x_0 =(d_1 d_2 d_3 -c_1 c_2 c_3)$$
satisfies $x_0 < c_1$, then $d_1 < c_1$. 
\end{proposition} 

\it Proof. \rm \ Let $P(x)=x^3 -\alpha x^2 +(c_1 c_2 +c_1 c_3 +c_2 c_3) x -c_1 c_2 c_3$ and $\tilde{P}(x)=
x^3 -\alpha x^2 +(d_1 d_2 +d_1 d_3 +d_2 d_3) x -d_1 d_2 d_3$. The function 
  $$Q(x)=\tilde{P}(x)-P(x)=\Big((d_1 d_2 +d_1 d_3 +d_2 d_3) -(c_1 c_2 +c_1 c_3 +c_2 c_3)\Big) x +(c_1 c_2 c_3 -d_1 d_2 d_3)$$
has the unique zero point $x_0$. The function $Q(x)$ is strictly positive for $x > x_0$ and hence $\tilde{P}(x) > P(x)$ on
 the half-line $x_0 <x <\infty$.  Thus assumption  $x_0 < c_1$ implies that $d_1 < c_1$ by Lemma 2.4(b). \qed

\medskip
By Proposition \ref{cubic3}, we obtain the following.

 \begin{corollary} \label{cubiccriterion}
Let $\bc=(c_1, c_2, c_3)$ and $\bd=(d_1, d_2, d_3)$ be positive vectors with descending coordinates such that 
  $\alpha=c_1 +c_2 +c_3 =d_1+ d_2+ d_3$ and $c_1 c_2 c_3 \ne d_1 d_2 d_3$.
If $c_1 c_2 +c_1 c_3 +c_2 c_3 < d_1 d_2 +d_1 d_3 +d_2 d_3$ and the value $x_0$ determined by 
  $$\Big((d_1 d_2 +d_1 d_3 +d_2 d_3) -(c_1 c_2 +c_1 c_3 +c_2 c_3)\Big) x_0 =(d_1 d_2 d_3 -c_1 c_2 c_3)$$
satisfies
  $$x_0 \leq \frac{1}{3}(c_1 +c_2 +c_3 +\sqrt{c_1^2 +c_2^2 +c_3^2 -c_1 c_2-c_1 c_3 -c_2 c_3}),$$
then $d_1 < c_1$. Therefore if $c_1 c_2 +c_1 c_3 +c_2 c_3 < d_1 d_2 +d_1 d_3 +d_2 d_3$ and 
   $$x_0 \leq \frac{1}{3}(c_1 +c_2 +c_3),$$
then  $d_1 < c_1$. 
 \end{corollary}

{\it Proof.} The cubic polynomial 
   $$P(x)=(x-c_1)(x-c_2)(x-c_3)=x^3 -(c_1 +c_2 +c_3) x^2 +(c_1 c_2 +c_1 c_3 +c_2 c_3) -c_1 c_2 c_3$$
has the derivative 
  $$P'(x)=3x^2 -2(c_1 +c_2 +c_3) x +(c_1 c_2 +c_1 c_3 +c_2 c_3).$$
 The quadratic equation $P'(x)=0$ has the two real roots 
   $$Z(\epsilon)=\frac{1}{3}(c_1 +c_2 +c_3 +
   \epsilon (c_1^2 +c_2^2 +c_3^2 -c_1 c_2 -c_1 c_3 -c_2 c_3)^{1/2}),$$
where $\epsilon=\pm 1$. Clearly, we have 
   $$0 <c_3 < Z(-1) < c_2 <Z(+1) < c_1.$$
Hence $x_0 \leq Z(+1)$, and thus $x_0 < c_1$. Then conclusion $d_1<c_1$ follows from Proposition \ref{cubic3}.
\qed

\medskip

Let $0 <a_1 <a_2 <a_3 <a_4 < a_5 < a_6$. For $j=1,2,\ldots, 60$, $j\ne 21$,  assume that $\eta\in S_6$ is the corresponding permutation which gives rise to the cyclic shift matrix $S(j)$, i.e., 
$S(j)=S(a_{\eta(1)}, \ldots, a_{\eta(6)})$.  Define the polynomial
    $$G_j(x)=\det(xI_6 -2\Re(e^{i\theta} S(a_{\eta(1)}, a_{\eta(2)}, a_{\eta(3)}, a_{\eta(4)}, a_{\eta(5)}, a_{\eta(6)})))$$
    $$ - \det(x I_6 -2\Re(e^{i \theta} S(21))).$$
Then  $G_j(x)=\alpha(j) x^2 +\beta(j)$, where 
 $\alpha(j)=S(21)-S(j)$, the difference of cyclic sums,  and 
 $\beta(j)=a_1^2 a_4^2 a_5^2 +a_2^2 a_3^2 a_6^2 -a_{\eta(1)}^2 a_{\eta(3)}^2 a_{\eta(5)}^2 -a_{\eta(2)}^2 a_{\eta(4)}^2 a_{\eta(6)}^2$. 
According to Lemma \ref{top21}, $S(21)>S(j)$ for all $j$, and thus $\alpha>0$. Denote $X(j)=-\beta(j)/\alpha(j)$. Then $G_j(x)>0$ whenever $x^2>X(j)$.

\medskip
The following result is useful for the comparison of the numerical range of
$S(21)$ and  the numerical ranges of $S(j)$ for 
$j=9,13,15,19,20,33,37,39,45$.  

\medskip

\begin{lemma}\label{n=6intersection}
 Let $0 < a_1 <a_2 <a_3 <a_4 <a_5 < a_6$. 
Then for $j=9, 13, 15, 19,$ $ 20, 33, 37, 39, 45$, 
$$X(j) \leq X_0 +\frac{r_5}{3}=\frac{1}{3}(a_1^2 +a_2^2 +a_3^2 +2a_5^2 +a_6^2) <\Big(\lambda_1(2\Re(e^{i \theta}S(21)))\Big)^2,$$ 
$0 \leq \theta \leq 2\pi$, 
where $X_0=\frac{1}{3}(a_1^2 +a_2^2 +a_3^2 +a_4^2+a_5^2 +a_6^2)$.
\end{lemma} 

\it Proof. \rm \ First, the values $X(13), X(33), X(37), X(39), X(45)$ are computed and respectively given by 

$\frac{1}{r_4}\Big(r_1 r_4 -r_2 (r_1 +r_2+r_3)\Big)$,\   

$\frac{1}{r_2 +r_3 +r_5 +r_6} 
\Big(r_1 r_2 + r_2^2 + r_1 r_3 + 2 r_2 r_3 + r_3^2 + r_2 r_4 + r_3 r_4 +(r_1+ r_2 + r_3)(r_5+ r_6)\Big)$,\   

$\frac{r_5 +r_6}{r_4 r_6 +r_5(r_3+r_4+r_5+r_6)}
\Big(r_1(r_4 +r_5) -r_2(r_1 +r_2 +r_3)\Big)$,\ 

$\frac{r_4+r_5}{r_4(r_2 +r_6) +r_5(r_2+r_3+r_4+r_5 +r_6)}
\Big(r_1r_2 + r_2^2 + r_2 r_3 + r_2 r_4 + r_1 r_5 + r_2 r_5 + r_1 r_6 + r_2 r_6\Big)$,\ 

$\frac{1}{r_3 +2r_4 +r_5}\Big(r_1 r_3 + r_2 r_3 + r_3^2 + 2 r_1 r_4 + 2 r_2 r_4 + 2 r_3 r_4 + 
r_4^2 + r_1 r_5 +  r_2 r_5 + r_3 r_5 + r_4 r_5\Big)$. 

\medskip\noindent
Moreover, we see that all of the terms computed satisfy  
 $X(j) \leq X_0$. 

\medskip
 Second, the values $X(19), X(15), X(20), X(9)$ are respectively computed by  

\medskip 
 $\frac{1}{r_4}\Big(r_1 r_4 + r_2 r_4 + r_3 r_4 + r_4^2 + r_4 r_5 - r_1 r_6 - r_2 r_6 - r_3 r_6\Big)$, \  
 
 $\frac{1}{r_2 +r_6}\Big(r_1 r_6 +r_2 (r_1+r_2+r_3+r_4+r_5+r_6)\Big)$,\ 
 
 $\frac{r_2 +r_3}{r_2 r_4 +r_3(r_2 +r_3+r_4 +r_5)}
 \Big(r_1 r_3 + r_2 r_3 + r_3^2 +
 r_1 r_4 + r_2 r_4 + 2r_3 r_4 + r_4^2 + r_3 r_5 +r_4 r_5 - r_1 r_6 - r_2 r_6\Big)$,\ 
 
 $\frac{r_3+r_4}{r_4 (r_2 +r_6) +r_3 (r_2+r_3+r_4+r_5+r_6)}\Big(r_1(r_2+r_3+r_6)+r_2
(r_2+2r_3)+r_3^2 +(r_2+r_3)(r_4+r_5+r_6)\Big)$.

\medskip\noindent
To evaluate these values, we set  
 \begin{eqnarray*}
&& \tilde{X}(19)=X_0 +\frac{r_5}{3}=\frac{1}{3}(6r_1+5r_2+4r_3+3r_4+2r_5+r_6) +\frac{r_5}{3},\\
&& \tilde{X}(15) =X_0 +\frac{r_2 r_5}{3(r_2 +r_6)},\\
&& \tilde{X}(20)=X_0 +\frac{r_2 r_4 r_5}{3(r_2 r_4 +r_2 r_3+r_3^2+r_3 r_4+r_3 r_5)},\\
&& \tilde{X}(9)=X_0 +\frac{r_2 r_4 r_5}{3(r_2 r_4 +r_2 r_3 +r_3^2 +r_3 r_4 +r_3 r_5 +r_3 r_6 +r_4r_6)}.
 \end{eqnarray*}
Then $\tilde{X}(j) \leq X_0 +r_5/3$ for $j=19,15,20,9$.

\medskip
We clearly have $X(19) \leq \tilde{X}(19)$. Furthermore,
  \begin{eqnarray*}
  && 3(r_2 +r_6)\Big(\tilde{X}(15) -X(15)\Big)=r_2 (3r_1 +2r_2+r_3) +r_6(3r_1+3r_2+4r_3+2r_5+r_6) \geq 0,\\
  && 3(r_2 r_3 +r_3^2 +r_2 r_4 +r_3 r_4 +r_3 r_5)\Big(\tilde{X}(20)-X(20)\Big)=(r_3+r_4) (r_2 +r_3)(3r_1 +2r_2 +r_3) \\
  && \hskip20pt +r_3 r_5 (6r_1+4r_2+3r_3 +2r_4 +2r_5 +r_6) +r_6 (r_2 +r_3)(3r_1 +3r_2 +r_3+r_4) \geq 0,\\
  && 3(r_2 r_3 +r_3^2 +r_2 r_4 +r_3 r_4 +r_3 r_5 +r_3 r_6 +r_4 r_6) \Big(\tilde{X}(9) -X(9)\Big)=(r_2 +r_3)(r_3 +r_4)(3r_1 +2r_2 +r_3)\\
  &&  \hskip20pt +r_3 r_5 (6r_1 +4r_2 +3r_3 +2r_4 +2r_5 +3r_6)+r_6 (3r_1 r_3+3r_2 r_3 +2r_3^2 +3r_1 r_4  \\
  &&  \hskip20pt +3r_2 r_4 +5r_3 r_4 +3r_4^2 + 2r_4 r_5 +r_3 r_6 +r_4 r_6) \geq 0.
   \end{eqnarray*}
It follows that
  $$X(j) \leq \tilde{X}(j) \leq X_0 +{r_5\over 3},$$
for $j=15,20,9$.  This proves that  
$X(j) \leq X_0 +\frac{r_5}{3}$ for all $j=9, 13, 15, 19,$ $ 20, 33, 37, 39, 45$.  

\medskip
The cyclic shift matrix $S(21)=S(a_4,a_2, a_1, a_3,a_5, a_6)$ is equivalent to 
 $S(a_2,a_4,a_6,a_5,a_3, a_1)$. By the fact that the numerical range of a principal 
 submatrix is contained in the numerical range of the matrix, we have the inclusion:
    $$W(S(a_4,a_6,a_5,0)) \subseteq  
    W(S(a_4,a_6,a_5,a_3,0)) \subseteq W(S(a_2,a_4,a_6,a_5,a_3,a_1)).$$
 We compute that 
 $${\rm det}(x I_4  -2\Re(S(a_4,a_6,a_5,0)))=x^4 -(a_4^2 +a_5^2 +a_6^2)x^2 +a_4^2 a_5^2.$$
Hence
 \begin{eqnarray*}
&&\Big( \lambda_1(2 \Re(e^{i\pi/6} S(a_2,a_4,a_6,a_5,a_3,a_1)))\Big)^2\\
&& \ge \frac{1}{2}\Big(a_4^2 +a_5^2 +a_6^2 +\sqrt{(a_4^2 +a_5^2 +a_6^2)^2 -4a_4^2 a_5^2}\Big)\\
&& =\frac{1}{2}\Big(a_4^2 +a_5^2 +a_6^2 +\sqrt{a_4^4 +a_5^4 +a_6^4 -2a_4^2 a_5^2 +2a_4^2 a_6^2 +2a_5^2 a_6^2}\Big)\\
&&  \ge \frac{1}{2}\Big(a_4^2 +a_5^2 +a_6^2  +\sqrt{a_4^4 +a_5^4 +a_6^4 -2a_4^2 a_5^2 -2a_4^2 a_6^2 +2a_5^2 a_6^2}\Big)\\
&&=a_5^2 +a_6^2\\
&&\ge   \frac{1}{3}(a_5^2 +a_6^2 +a_6^2 +2a_5^2 +a_6^2)\\ 
&& > \frac{1}{3}(a_1^2 +a_2^2 +a_3^2 +2a_5^2 +a_6^2).
\end{eqnarray*}
\qed

\medskip 

Now, we are ready to present the following. 

\medskip\noindent
{\bf Proof of Theorem \ref{upper=bound45678}}
For $j\in \{9, 13,15,19,20,33,37,39,45\}$,  assume that $\eta\in S_6$ is the corresponding permutation which gives rise to the cyclic shift matrix $S(j)$.
Given $\theta \in [0, 2\pi]$, 
 $${\rm det}(xI_6 -2\Re(e^{i \theta}S(j))) -{\rm det}(xI_6 -2\Re(e^{i \theta}S(21)))$$
 $$=(S(21) -S(j)) x^2 +a_1^2 a_3^2 a_5^2 +a_2^2 a_4^2 a_6^2  -a_{\eta(1)}^2 a_{\eta(3)}^2 a_{\eta(5)}^2 -a_{\eta(2)}^2 a_{\eta(4)}^2 a_{\eta(6)}^2>0$$
for $x^2> X(j)$. 
 By Lemma \ref{n=6intersection},
 $X(j) < \Big(\lambda_1(2 \Re(e^{i \theta}S(21))\Big)^2$. Then, according to  Corollary \ref{cubiccriterion}, we obtain that 
    $$ \lambda_1(2\Re(e^{i \theta} S(j)) \leq \lambda_1(2 \Re(e^{i \theta} S(21)),$$
and thus $W(S(j)) \subseteq W(S(21))$ for $j=9,13, 15, 19, 20,33, 37, 39, 45$. Combining this 
inclusion with Proposition \ref{reduction}, we conclude that $W(S(j)) \subseteq W(S(21))$ for $j=1, \ldots, 60$. 
\qed

\section{Related results and further research}

Note that the fourth author has  verified Conjecture \ref{conj1}  for the cases $n= 7,8$, by 
extending  the proof for the case $n = 6$ with the help of the symbolic calculation program  Mathematica.
The computation is very involved.  It would be interesting to prove or disprove the
conjecture for general $n$. Of course, a related question is to determine the 
permutation $\tau$ so that $A_\tau$ yields the smallest 
$W(A_\tau)$ for $n \ge 6$, if such a $\tau$ exists.
The fourth author also has some partial results in this direction for the case that $n = 6$.

As pointed out in Section 1, our study is closely related to the maximum and
minimum eigenvalues of the adjacency matrix of a given graph. In our case, 
 if $A = S(a_1, \dots, a_n)$ with positive numbers $a_1, \dots, a_n$,  then $A+A^t$ is the adjacency matrix of a weighted undirected cycle. One may consider a bound  
$\gamma$ for the largest
eigenvalue of $A+A^t$ for a general nonnegative matrix $A$ and interpret it in terms of the 
adjacency matrix of a weighted graph; then the numerical radius of $A$ will satisfy  $r(A) \le \gamma$.  
For example, that problem is solved in the case that $A+A^t$ corresponds to the adjacency 
matrix of a tree (and therefore a forest); see \cite{Cheung}. 
Of course, in general it is more challenging to determine whether
$W(A_\sigma)$ is maximum or minimum in terms of set inclusion.

Our study is also related to the study of positive linear maps on operator 
algebras and operator subspaces (systems).
Recall that a linear map $\Phi: V_1 \rightarrow V_2$ from a subspace $V_1$ of $M_n$ to
a subspace $V_2$ of $M_m$ is $k$-positive for $k \in \IN$ if
$(\Phi(A_{ij}))\in M_n(V_2)$ is positive semidefinite whenever 
$(A_{ij}) \in M_n(V_1)$ is positive semidefinite;
$\Phi$ is completely positive if $\Phi$ is $k$-positive for all $k \in \IN$. 
One may connect the study of the inclusion of the numerical range 
with positive linear maps and completely positive maps; see \cite{CL}.

\begin{proposition} 
Let $A, B \in M_n$. Then $W(B) \subseteq W(A)$ if and only if any one of the following holds.
\begin{itemize}
\item [{\rm (a)}]
$\mu I - \Re(e^{i\theta }B)$ is positive definite whenever
$\mu I - \Re(e^{i\theta }A) $ is positive definite.
\item [{\rm (b)}] The map $T: \span \{I, A, A^*\} \rightarrow \span\{I, B, B^*\}$ 
defined by $\alpha I + \beta A + \gamma A^* \mapsto \alpha I + \beta B + \gamma B^*$
send positive definite matrices to positive semidefinite matrices.
\item [{\rm (c)}] 
$\lambda_1(\Re(e^{i\theta }A)) \ge \lambda_1(\Re(e^{i\theta }B))$ for all $\theta \in [0, 2\pi)$.
\end{itemize}
In the case that $A, B$ are Hermitian, the above conditions are  equivalent to the following.
\begin{itemize}
\item[{\rm (d)}] There is a unital completely positive linear map $\Psi: \span\{I, A\}
\rightarrow \span\{I, B\}$ such that 
$\Psi(A) = B$.
\end{itemize}
\end{proposition}

Consequently, our results on $W(B) \subseteq W(A)$ for 
$B = S(b_1, \dots,b_n)$ and $A = S(a_1, \dots, a_n)$ provide
necessary and sufficient conditions for the existence of a positive 
linear map from $\{I, A, A^*\}$ to $\{I, B, B^*\}$. Note that
there is a positive linear map from $\{I, A, A^*\}$ to $\{I, B, B^*\}$
if and only if there is a positive linear map from 
$\{I, U^*AU, U^*A^*U\}$ to 
$\{I, V^*BV, V^*B^*V\}$ for any unitary matrices $U$ and $V$.
One may also consider matrices of different sizes, or infinite sizes
if we replace $W(X)$ by its closure. So, there might be interest in 
extending our results to general operators and infinite graphs.

\begin{remark} \rm One might seek to extend Conjecture \ref{conj1} to find the minimal numerical range. Suppose that $n \ge 2,$  and  $0 \le a_1 \le \cdots \le a_n$. The reviewer of this paper proposes that
\begin{equation}\label{conjecture2}
W(B) \subseteq W(A_\sigma) \qquad \hbox{ for any permutation } \ \sigma,
\end{equation}
where  $B = S(\dots,a_{n-4},a_4,a_{n-2},a_2,a_n, a_1,a_{n-1},a_3,a_{n-3},a_5,a_{n-5},\dots).$

\medskip

The inclusion \eqref{conjecture2} is confirmed for $n=4$ and $n=5$ in Theorem \ref{n=4} and Theorem \ref{main-5} respectively. 
We provide an example to show that the inclusion \eqref{conjecture2} is not true for $n=6$. 

By letting $a_1=0, a_2^2=3, a_3^2=4, a_4^2=8, a_5^2=13, a_6^2=30$, we compare the numerical ranges $W(S(a_4, a_2,a_6, a_1,a_5,a_3))$ and $W(S(a_3,a_2,a_6,a_1,a_5,a_4))$. Easily, we compute that
\begin{equation}\label{ce1} \det(x I_6 -(S(a_4, a_2,a_6, a_1,a_5,a_3) +S(a_4, a_2,a_6, a_1,a_5,a_3)^T))=x^6 -58 x^4  +905 x^2 -3120.\end{equation}
On the other hand,
\begin{equation}\label{ce2}\det(x I_6 -(S(a_3, a_2,a_6, a_1,a_5,a_4) +S(a_3, a_2,a_6, a_1,a_5,a_4)^T))=x^6 -58 x^4  +865 x^2 -1560.\end{equation}
The largest zeros of the polynomials (\ref{ce1}) and (\ref{ce2}) are, numerically approximated by 5.8493 and 5.8067 respectively. This fact disproves the conclusion
 $$W(S(a_4, a_2,a_6, a_1,a_5,a_3))\subseteq W(S(a_3,a_2,a_6,a_1,a_5,a_4)).$$
\end{remark}

\section*{Acknowledgments}
The authors would like to express their thanks to an anonymous referee for his/her careful reading and valuable comments which  led to the improvement of the present version of the paper, in particular, for proposing the Lagrange method for proving Theorem \ref{3.2} and simplifying the proof of Proposition \ref{prop3}.
Kirkland's research is supported in part by NSERC
Discovery Grant RGPIN-2019-05408. 
Li is an affiliate member of the Institute for Quantum Computing, University of Waterloo;
his research was partially supported by the Simons Foundation Grant 851334.

\bigskip

\noindent
{\bf Declaration of competing interest}

The authors declare no conflict of interest.

\end{document}